\documentclass[11pt,leqno]{article}

\usepackage{amsmath,amsfonts,amscd,amssymb,theorem}

\long\def\comment#1\endcomment{}

\makeatletter
\begingroup
\gdef\th@dotted{\normalfont\itshape
  \def\@begintheorem##1##2{%
        \item[\hskip\labelsep \theorem@headerfont ##1\ ##2.]}%
\def\@opargbegintheorem##1##2##3{%
   \item[\hskip\labelsep \theorem@headerfont ##1\ ##2\ (##3).]}}
\endgroup
\makeatother

\theoremstyle{dotted}

\newtheorem{theorem}{Theorem}[section]
\newtheorem{lemma}[theorem]{Lemma}

\newtheorem{prop}[theorem]{Proposition}
\newtheorem{corr}[theorem]{Corollary}

\makeatletter
\begingroup
\gdef\th@upshape{\normalfont
  \def\@begintheorem##1##2{%
        \item[\hskip\labelsep \theorem@headerfont ##1\ ##2.]}%
\def\@opargbegintheorem##1##2##3{%
   \item[\hskip\labelsep \theorem@headerfont ##1\ ##2\ (##3).]}}
\endgroup
\makeatother

\theoremstyle{upshape}

\newtheorem{defn}[theorem]{Definition}
\newtheorem{remark}[theorem]{Remark}
\newtheorem{exa}[theorem]{Example}

\makeatletter
\renewcommand{\subsection}{\@startsection{subsection}{2}{0pt}{-3ex
plus -1ex minus -0.2ex}{-2mm plus -0pt minus
-2pt}{\normalfont\bfseries}} 
\renewcommand{\subsubsection}{\@startsection{subsubsection}{3}{0pt}{-3ex
plus -1ex minus -0.2ex}{-2mm plus -0pt minus
-2pt}{\normalfont\bfseries}} 
\makeatother

\makeatletter
\@addtoreset{equation}{section}
\makeatother

\newcommand{\cntrct}                
{\hspace{2pt}\raisebox{1pt}{\text{$\lrcorner$}}\hspace{2pt}}

\newcommand{\proof}[1][Proof.]{\smallskip\noindent{\em #1}}
\def\endproof{\hfill\ensuremath{\square}\par\medskip}

\renewcommand{\labelenumi}{{\normalfont(\roman{enumi})}}

\def\eqref#1{\thetag{\ref{#1}}}

\let\latexref=\ref
\def\ref#1{{\normalfont{\latexref{#1}}}}

\newcommand{\wt}{\widetilde}

\newcommand{\dg}{\dagger}

\newcommand{\idot}{{\:\raisebox{1pt}{\text{\circle*{1.5}}}}}
\newcommand{\eps}{\varepsilon}
\newcommand{\ups}{\upsilon}
\renewcommand{\phi}{\varphi}

\newcommand{\Hom}{\operatorname{Hom}}

\newcommand{\Fun}{\operatorname{Fun}}

\newcommand{\Sec}{\operatorname{Sec}}
\newcommand{\Ch}{\operatorname{Ch}}

\newcommand{\id}{\operatorname{\sf id}}

\newcommand{\A}{\mathcal{A}}
\newcommand{\C}{\mathcal{C}}

\newcommand{\hash}{\sharp}

\newcommand{\Sets}{\operatorname{Sets}}
\newcommand{\Cat}{\operatorname{Cat}}

\newcommand{\ppt}{{\sf pt}}

\newcommand{\Ind}{\operatorname{Ind}}
\newcommand{\Comp}{\operatorname{Comp}}

\newcommand{\copr}{\sqcup}

\newcommand{\N}{{\mathbb N}}

\newcommand{\E}{\mathcal{E}}

\newcommand{\colim}{\operatorname{\sf colim}}
\renewcommand{\lim}{\operatorname{\sf lim}}

\newcommand{\V}{{\sf V}}

\newcommand{\Y}{{\sf Y}}
\newcommand{\YY}{\mathcal{Y}}
\newcommand{\X}{\mathcal{X}}
\newcommand{\Z}{\mathcal{Z}}

\newcommand{\Cone}{\operatorname{Cone}}

\newcommand{\ev}{\operatorname{\sf ev}}

\newcommand{\Pos}{\operatorname{Pos}}

\newcommand{\Tw}{\operatorname{\sf tw}}

\newcommand{\hht}{\operatorname{ht}}

\newcommand{\Filt}{\operatorname{Filt}}

\newcommand{\trl}{\operatorname{\triangleleft}}

\newcommand{\ssetminus}{\smallsetminus}

\newcommand{\ra}{\overrightarrow}
\newcommand{\la}{\overleftarrow}

\newcommand{\vtimes}{\ra{\times}}

\title{Taming large categories}

\author{D. Kaledin\thanks{Partially supported by the HSE University
    Basic Research Program}}

\date{{\em To Sasha Anan'in}}

\begin{document}

\maketitle

\tableofcontents

\section*{Introduction.}

In standard category theory, there are two types of categories:
small ones, whose objects and morphisms form a set, and large ones
that are allowed to have a proper class of objects. However, even
for a large category $\C$, morphisms between any two given objects
$c$ and $c'$ are always assumed to form a set. This has an
unpleasant side effect: functors between two large categories $\C$,
$\E$ do not necessarily form a well-defined category. Indeed, for
any two functors $E,E':\C \to \E$, morphisms from $E$ to $E'$ might
form a proper class.

The standard solution to this difficulty is to only consider
functors that preserve colimits of some sort, and use this fact to
insure that morphisms $E \to E'$ do form a set. However, in this
approach, one has to require the existence of these colimits, and
one obviously does not want to require too much. In particular, for
small source categories, there is no problem, so ideally, for a
small category, one does not want to require anything at all.

A general formalism that satisfies these requirement -- or rather,
almost satisfies them, up to a relatively mild operation of taking
Karoubi closure -- is given by the theory of accessible categories
and filtered colimits. This has been initiated by Grothendieck
\cite{SGA4} and developed into a fully fledged theory by Gabriel and
Ulmer \cite{GU}. After that, the theory lay dormant for a couple of
decades, and then was given its present form in \cite{MP}, and then
in \cite{AR}.

The present paper gives a short and self-contained overview of the
subject. There are no new results whatsoever, and we do not claim
much novelty in the proofs either. There are two features of our
exposition that are worth mentioning.
\begin{enumerate}
\item We systematically use ``categories of elements'' for functors
  to $\Sets$. The fact that this makes many proofs simpler and
  easier to comprehend is well-known in folklore, but in our
  opinion, somewhat missing in the available literature. For reasons
  of brevity, we have avoided the general machinery of Grothendieck
  fibrations of \cite{SGA}, although morally, categories of elements
  are just discrete Grothendieck fibrations, and should be treated
  as such.
\item There are good reasons to treat the theory of accessible
  categories as part of Logic; this point of view is very prominent
  in \cite{MP}, and present to some extent in \cite{AR}. However,
  there are also good reasons {\em not} to do this, and this is the
  approach adopted here. We work in a purely categorical context,
  such as e.g.\ \cite{SGA4}, and do not mention Logic at all.
\end{enumerate}
Having said all that, let us emphasize that especially \cite{AR} is
a beautiful book that gives a comprehensive and thorough exposition
of the subject, well-written and reasonably concise, so one should
really justify the necessity of one more overview. Our original
motivation for doing this was our current work-in-progress on
homotopical enhancements for category theory based on Grothendieck's
idea of a ``derivator'', \cite{big}. In that context, the machinery
of accessible categories is even more indispensable than in the
usual category theory, since without it, even fibered products are
not well-defined. Of course, the whole theory has to be re-done in
the enhanced setting, but this turns out to be suprisingly
straightforward -- with one or two technical exceptions, all of the
arguments presented here extend to the enhanced context almost {\em
  verbatim}. Our second motivation is the sheer beauty and
simplicity of the theory. Modulo general categorical preliminaries,
quite standard and needed mostly to fix notation and terminology,
everything fits together nicely in about 15 pages. Once we have
realized that this can be done, it was hard to resist the temptation
to actually do it. To keep things really short, we completely avoid
the theory of presentable categories that usually comes packaged
together with the material here, and concentrate entirely on the
accessible case.

\subsection*{Acknowledgements.} Both in preparation of this
particular paper, and in thinking about category theory in general,
I have benefited from discussions with many colleagues. I am
especially grateful to A. Efimov and D. Tereshkin, and also to
Yu. Berest, M. Kapranov and A. Polischuk. I am grateful to the
referee for their careful reading of the manuscript and valuable
comments and suggestions.

\bigskip

\noindent
This paper is dedicated to Sasha Anan'in.

\bigskip

I never knew Sasha in the ancient days before the demise of the
USSR; he was in Novosibirsk and then left for Brazil, while I was in
Moscow and then went to the US for my graduate studies. He
re-appeared in the horizon somewhere around 2010, out of the blue,
and immediately became quite a presence in Moscow, both virtually,
throughout all the subsequent years, and physically, during one
short year he actually spent at our Laboratory of Algebraic Geometry
at HSE University. He is an honorary member of the Laboratory, for
as long as it lasts. He is a dear friend, who is sorely missed. His
enthusiasm for mathematics -- all kinds of mathematics, since it is
all connected, when one knows how to look -- was absolutely
infectuous, and hugely encouraging. This includes everything related
to our craft -- doing mathematics, teaching mathematics, talking
about mathematics, simply enjoying it, in English, Russian, Spanish or
Portuguese. This paper is a chapter of a general textbook on
category theory that we hoped to write together, at some time,
somewhere, in some other life. !`Hasta siempre!

\section{Preliminaries.}\label{cat.sec}

\subsection{Categories and functors.}

We work in the minimal set-theoretic setup limited to small and
large categories; even for a large category, $\Hom$-sets are
small. We use ``map'', ``morphism'' and ``arrow''
interchangeably. For any category $\C$, we write $c \in \C$ as
shorthand for ``$c$ is an object of $\C$''. For any category $\C$,
we denote by $\C^o$ the opposite category, and for any functor
$\gamma:\C_0 \to \C_1$, we denote by $\gamma^o:\C_0^o \to \C_1^o$
the opposite functor. We denote by $\ppt$ the point category (a
single object, a single morphism), and for any $c \in \C$, we denote
by $\eps(c):\ppt \to \C$ the embedding onto $c$.

A functor $\gamma:\C \to \C'$ is {\em full} resp.\ {\em faithful} if
it surjective resp.\ injective on $\Hom$-sets, and {\em essentially
  surjective} if it is surjective on isomorphism classes of
objects. A functor $\gamma$ is {\em conservative} if any map $f$
with invertible $\gamma(f)$ is itself invertible. A functor is an
{\em equivalence} if it is fully faithful and essentially
surjective, or equivalently, if it is invertible up to an
isomorphism. A category is {\em essentially small} if it is
equivalent to a small category. We assume known the notion of a
right or left-adjoint functor, and the fact that adjoints are unique
(if they exist). In particular, an object $o \in \C$ is {\em
  initial} resp.\ {\em terminal} iff $\eps(o):\ppt \to \C$ is left
resp.\ right-adjoint to the tautological projection $\C \to \ppt$,
and an initial resp.\ terminal object is unique, if it exists.

A {\em cardinal} is an isomorphism class of sets, and for any set
$S$, we denote the corresponding cardinal by $|S|$. Cardinals are
ordered in the usual way ($|S| \leq |S'|$ iff there exists an
injective map $S \to S'$). The class of all cardinals is
well-ordered, and for any cardinal $\kappa$, the {\em successor
  cardinal} $\kappa^+$ is the smallest cardinal $\kappa^+ >
\kappa$. A cardinal $\kappa$ is {\em regular} if for any map $f:S
\to S'$ of sets such that $|S'| < \kappa$ and $|f^{-1}(s)| < \kappa$
for any $s \in S'$, we have $|S| < \kappa$. All regular cardinals
are infinite, and for any infinite cardinal $\kappa$, the successor
cardinal $\kappa^+$ is regular.

\begin{exa}\label{reg.exa}
Not all infinite cardinals are regular. For example, starting from
some infinite cardinal $\kappa = \kappa_0$, we can define $\kappa_n
= \kappa_{n-1}^+$, $n \geq 0$, and then $\kappa_\infty = \cup_{n
  \geq 0}\kappa_n$ is not regular.
\end{exa}

For any $n \geq 0$, a {\em chain} of length $n$ in a category $\C$
is a diagram
\begin{equation}\label{c.idot}
\begin{CD}
c_0 @>>> \dots @>>> c_n,
\end{CD}
\end{equation}
and a chain is {\em non-degenerate} if none of the maps $c_{l-1} \to
c_l$, $1 \leq l \leq n$ in \eqref{c.idot} is an identity map. If
$\C$ is small, then all non-degenerate chains in $\C$ form a set
$\Ch(\C)$, and we denote $|\C| = |\Ch(\C)|$. In particular, chains
of length $0$ are objects in $\C$, and chains of length $1$ are
morphisms. If the subset $\Ch_{\leq 1}(\C) \subset \Ch(\C)$ of
chains of length $\leq 1$ is infinite, then $|\C| = |\Ch_{\leq
  1}(\C)|$, but in general, this need not be true.

\begin{exa}\label{P.exa}
Let $P$ be the category with one object and one non-identity
morphism $p$ such that $p^2=p$. Then $P$ has exactly one
non-degenerate chain of each length $n$, so $|P|$ is infinite, while
$|\Ch_{\leq 1}(P)| = 2$.
\end{exa}

A small category is {\em discrete} if all its maps are identity
maps, so that a discrete small category is the same thing as a
set. We denote by $\Sets$ resp.\ $\Cat$ the categories of sets
resp.\ small categories. A category $\C$ is {\em connected} if it is
not empty, and does not admit a decomposition $\C \cong \C_1 \copr
\C_2$ with non-empty $\C_1$, $\C_2$ (or equivalently, any two
objects $c,c' \in \C$ can be connected by a finite zigzag of
morphisms). We have the fully faithful embedding $\Sets \subset
\Cat$ identifying discrete categories and sets, and it has a
left-adjoint functor $\pi_0:\Cat \to \Sets$ sending a small category
to the set of its connected components. If a functor $\gamma$ admits
an adjoint, then $\pi_0(\gamma)$ is an isomorphism. For any regular
cardinal $\kappa$, we say that a set or a small category $X$ is {\em
  $\kappa$-bounded} if $|X| < \kappa$, and we denote by
$\Sets_\kappa \subset \Sets$ resp.\ $\Cat_\kappa \subset \Cat$ the
full subcategories spanned by $\kappa$-bounded sets resp.\ small
categories.

An object $c \in \C$ in a category $\C$ is a {\em retract} of an
object $c' \in \C$ if it is equipped with maps $a:c \to c'$, $b:c'
\to c$ such that $b \circ a = \id$. In such a situation, the
endomorphism $p=a \circ b:c' \to c'$ of the object $c'$ is
idempotent, that is, $p^2=p$. An idempotent endomorphism of an
object is a {\em projector}, and $c$ is the {\em image} of the
projector $p$. The image of a projector is unique, if it exists, and
it is automatically preserved by all functors $\C \to \C'$ to any
category $\C'$. A category $\C$ is {\em Karoubi-closed} if all the
projectors have images. The {\em Karoubi closure} $P(\C)$ of a
category $\C$ is the category of pairs $\langle c,p \rangle$, $c \in
\C$, $p:c \to c$ a projector, with maps $\langle c,p \rangle \to
\langle c',p' \rangle$ given by maps $f:c \to c'$ such that $p' \circ
f = f \circ p = f$. The category $P(\C)$ is Karoubi-closed, we have
the tautological full embedding $\eps:\C \to P(\C)$, $c \mapsto
\langle c,\id \rangle$, and it is universal: any functor $\C \to
\C'$ to a Karoubi-closed category $\C'$ factors through $\eps$,
uniquely up to a unique isomorphism. In the universal situation, the
Karoubi closure $P^=$ of the category $P$ of Example~\ref{P.exa} is
the category with two objects $o$, $o'$, with maps generated by maps
$a:o \to o'$, $b:o' \to o$, subject to relation $b \circ a =
\id$. The full subcategory in $P^=$ spanned by $o'$ is $P$, and $o$
is both the initial and the terminal object. An object in $P(\C)$ is
a functor $P \to \C$, and it lies in $\C \subset P(\C)$ iff the
functor extends to $P^=$.

For any categories $\C_0$, $\C_1$, $\C$ equipped with functors
$\gamma_l:\C_l \to \C$, $l=0,1$, the {\em lax fiber product} $\C_0
\vtimes_\C \C_1$ is the category of triples $\langle c_0,c_1,\alpha
\rangle$ of objects $c_l \in \C_l$, $l=0,1$ and a map
$\alpha:\gamma_0(c_0) \to \gamma_1(c_1)$, and the {\em product}
$\C_0 \times_\C \C_1 \subset \C_0 \vtimes_\C \C_1$ is the full
subcategory spanned by triples $\langle c_0,c_1,\alpha \rangle$ with
invertible $\alpha$. We denote by $\sigma:\C_0 \vtimes_\C \C_1 \to
\C_0$ resp.\ $\tau:\C_0 \vtimes_\C \C_1 \to \C_1$ the forgetful
functors sending $\langle c_0,c_1,\alpha \rangle$ to $c_0$
resp.\ $c_1$, and we write $\C_0 \times_\C \C_1 = \gamma_0^*\C_1$
when we want to emphasize its dependence on $\gamma_0$. A {\em
  commutative square} of categories and functors is a square
\begin{equation}\label{cat.sq}
\begin{CD}
\C_{01} @>{\gamma_{01}^1}>> \C_1\\
@V{\gamma_{01}^0}VV @VV{\gamma_1}V\\
\C_0 @>{\gamma_0}>> \C
\end{CD}
\end{equation}
of categories and functors equipped with an isomorphism $\gamma_1
\circ \gamma_{01}^1 \cong \gamma_0 \circ \gamma_{01}^0$. A square
\eqref{cat.sq} induces a functor $\C_{01} \to \C_0 \times_\C \C_1$,
and a square is {\em cartesian} if this functor is an
equivalence. For any functor $\pi:\C \to I$, the {\em left}
resp.\ {\em right comma-categories} are $\C /_\pi I = \C \vtimes_I I$
resp.\ $I \setminus_\pi \C = I \vtimes_I \C$, and we drop $\pi$ from
notation when it is clear from the context. The projection
$\sigma:\C / I \to \C$ has a fully-faithful left-adjoint $\eta:\C
\to \C/I$, $c \mapsto \langle c,\pi(c),\id \rangle$, and dually,
$\tau:I \setminus \C \to \C$ has a fully faithful right-adjoint
$\eta:\C \to I \setminus \C$; the functor $\pi$ itself then
decomposes as
\begin{equation}\label{comma.facto}
\begin{CD}
\C @>{\eta}>> \C /_\pi I @>{\tau}>> I,\qquad
\C @>{\eta}>> I \setminus_\pi \C @>{\sigma}>> I.
\end{CD}
\end{equation}
For any object $i \in I$, the {\em fiber} $\C_i$ is given by $\C_i =
\eps(i)^*\C$, and the {\em left} and {\em right comma-fibers} are
the fibers $\C /_\pi i = (\C /_\pi I)_i$, $i \setminus_\pi \C = (I
\setminus_\pi \C)_i$ of the projections $\tau$ resp.\ $\sigma$; the
projections $\sigma$ resp.\ $\tau$ then induce functors
$\sigma(i):\C / i \to \C$, $\tau(i):i \setminus \C \to
\C$. Explicitly, $\C / i$ resp.\ $i \setminus \C$ is the category of
pairs $\langle c,\alpha \rangle$, $c \in \C$, $\alpha$ a map $\pi(c)
\to i$ resp.\ $i \to \pi(c)$, and $\sigma(i)$ resp.\ $\tau(i)$ sends
such a pair to $c$.

\subsection{Categories of functors.}

For any category $\C$ and essentially small category $I$, functors
$I \to \C$ form a category that we denote $\Fun(I,\C)$, and we
simplify notation by writing $I^o\C = \Fun(I^o,\C)$. If $\C$ is
equipped with a functor $\pi:\C \to I$, we define the category
$\Sec(I,\C)$ of {\em sections} of the functor $\pi$ by the cartesian
square
\begin{equation}\label{sec.sq}
\begin{CD}
\Sec(I,\C) @>>> \Fun(I,\C)\\
@VVV @VVV\\
\ppt @>{\eps(\id)}>> \Fun(I,I),
\end{CD}
\end{equation}
where the bottom arrow is the embedding onto $\id:I \to I$. For any
essentially small categories $I$, $I'$ equipped with a functor
$\gamma:I' \to I$, we denote by $\gamma^*:\Fun(I,\C) \to
\Fun(I',\C)$ the pullback functor obtained by precomposition with
$\gamma$. In general, $\Fun(I,\C)$ comes equipped with the
evaluation pairing $\ev:\Fun(I,\C) \times I \to \C$, and for any
$I'$, a pairing $\gamma:I' \times I \to \C$ factors as
\begin{equation}\label{fun.facto}
\begin{CD}
I' \times I @>{\wt{\gamma} \times \id}>> \Fun(I,\C) \times I
@>{\ev}>> \C
\end{CD}
\end{equation}
for a functor $\wt{\gamma}:I' \to \Fun(I,\C)$, unique up to a unique
isomorphism. In particular, any essentially small category $I$ comes
equipped with a $\Hom$-pairing $\Hom_I(-,-):I^o \times I \to \Sets$,
and \eqref{fun.facto} then gives rise to the fully faithful {\em
  Yoneda embedding}
\begin{equation}\label{yo.eq}
\Y:I \to I^o\Sets.
\end{equation}
For any category $\C$, we denote by $\C^<$ resp.\ $\C^>$ the
category obtained by adding the new initial resp.\ terminal object
$o$ to $\C$, and we let $\eps:\C \to \C^>$, $\eps:\C \to \C^<$ be
the embeddings. For any two categories $\C_0$, $\C_1$, we let $\C_0
* \C_1 = (C_0^> \times \C_1^>) \ssetminus \{o \times o\}$, and we
note that we have
\begin{equation}\label{star.eq}
\C_0^> \times C_1^> \cong (C_0 * \C_1)^>.
\end{equation}
We also note that $\id:\C_0 \times \C_1 \to \C_0 \times \C_1$
extends to an embedding $(\C_0 \times \C_1)^> \to \C^>_0 \times
\C_1^>$, $o \mapsto o \times o$, and this has a left-adjoint
\begin{equation}\label{bicone.eq}
\C_0^> \times \C_1^> \to (\C_0 \times \C_1)^>.
\end{equation}
For any categories $I$, $\E$, a {\em cone} of a functor $E:I \to \E$
is defined as a functor $E_>:I^> \to \E$ equipped with an
isomorphism $\eps^*E_> \cong E$, and $E_>(o) \in \E$ is the {\em
  vertex} of the cone. Explicitly, a cone $E_>$ for $E$ is the same
thing as a pair $\langle e,\alpha \rangle$ of its vertex $e$ and a
morphism $E \to e$ to the constant functor $\C \to \E$ with value
$e$. For any such cone $E_> = \langle e,\alpha \rangle$, and any map
$f:e \to e'$ in $\E$, we denote $f_!E_> = \langle e',f \circ \alpha
\rangle$. All cones for a fixed $E$ form a category $\Cone(E)$. If
$\C$ is essentially small, then it can be described as
\begin{equation}\label{cone.eq}
\Cone(E) \cong E \setminus_{\gamma^*} \E,
\end{equation}
where $\gamma:I \to \ppt$ is the tautological projection. A cone is
{\em universal} if it is an initial object in $\Cone(E)$. If the
universal cone exists, then $E$ {\em has a colimit}, and the colimit
$\colim_I E$ is the vertex of this universal cone. If $I$ has a
terminal object $o \in I$, then $\colim_IE$ exists for any functor
$E:I \to \E$, and we have $\colim_IE \cong E(o)$. The Yoneda
embedding \eqref{yo.eq} has a colimit $\colim_I\Y \cong \ppt$ given
by the constant functor $\ppt:I^o \to \Sets$ with value $\ppt$. If
we have two categories $\E$, $\E'$, a functor $F:\E \to \E'$, and a
functor $E:I \to \E$ from an essentially small $I$ such that
$\colim_IE$ and $\colim_I(F \circ E)$ both exist, then the universal
property provides a map $\colim_I(F \circ E) \to F(\colim_IE)$, and
$F$ {\em preserves} the colimit $\colim_IE$ if this map is an
isomorphism (or equivalently, if $F \circ E_>$ for the universal
cone $E_>:I^> \to \E$ is universal).

\begin{exa}\label{P.cone.exa}
For the category $P$ of Example~\ref{P.exa}, with its Karoubi
closure $P^=$, the full embedding $\eps:P \to P^=$ admits a cone
$\eps_>:P^> \to P^=$ sending the terminal object $o$ to the terminal
object $o$.  This cone is universal and preserved by any functor
$P^= \to \E$ to any category $\E$, so that in particular, the image
of a projector in $\E$ is the colimit of the corresponding functor
$P \to \E$, and taking images of projectors commutes with
colimits. We also have the cone $a_!\eps_>:P^> \to P^=$ that
actually factors through $P$ and gives a cone for the identity
functor $\id:P \to P$; this cone is not universal.
\end{exa}

A category $\E$ is {\em cocomplete} if $\colim_IE$ exists for any
functor $E:I \to \E$ from a small $I$, and {\em $\kappa$-cocomplete}
for some regular cardinal $\kappa$ if this happens for any $I \in
\Cat_\kappa$ (we call such colimits {\em $\kappa$-bounded}). If $E$
is cocomplete or $\kappa$-cocomplete, then so is $\Fun(I,\E)$ for
any essentially small $I$.

\begin{exa}
The category $\Sets$ is cocomplete, and for any regular cardinal
$\kappa$, the category $\Sets_\kappa$ is $\kappa$-cocomplete (the
latter is essentially equivalent to the definition of a regular
cardinal).
\end{exa}

If we have a functor $\gamma:I \to I'$ between essentially small
categories, and a category $\E$ such that $\colim_{I / i}E$ exists
for any $i \in I'$ and $E:I/i \to \E$, then the pullback functor
$\gamma^*$ admits a left-adjoint {\em left Kan extension} functor
$\gamma_!$ given by
\begin{equation}\label{kan.eq}
f_!E(i) = \colim_{I/ i}\sigma(i)^*E, \qquad i \in I', E \in
\Fun(I,\E).
\end{equation}
In particular, this always happens if the target category $E$ is
cocomplete. If $I'=\ppt$, then \eqref{kan.eq} simply says that
$\colim_I$ is left-adjoint to the tautological functor $\E \to
\Fun(I,\E)$ sending $e \in \E$ to the constant functor $I \to \E$
with value $e$. The Yoneda embedding \eqref{yo.eq} can be described
in terms of left Kan extensions by $\Y(i) = \eps(i)^o_!\ppt$, $i \in
I$, and this shows that for any functor $\gamma:I_0 \to I_1$ between
essentially small categories, we have a commutative square
\begin{equation}\label{yo.sq}
\begin{CD}
I_0 @>{\Y}>> I_0^o\Sets\\
@V{\gamma}VV @VV{\gamma_!}V\\
I_1 @>{\Y}>> I_1^o\Sets.
\end{CD}
\end{equation}
Dually, the {\em limit} of a functor $E:I \to \E$ is $\lim_I E =
(\colim_{I^o}E^o)^o$, the {\em right Kan extension} $\gamma_*E$ with
respect to a functor $\gamma:I \to I'$ is given by $\gamma_*E =
(\gamma^o_!E^o)^o$, and both have the the dual versions of all the
properties of colimits and left Kan extensions. In particular, we
have the dual version of \eqref{kan.eq} for $\gamma_*$, with
$\colim_{I/i}$ replaced by $\lim_{i \setminus I}$. A category $\E$
is {\em complete} if $\E^o$ is cocomplete; this holds for
$\E=\Sets$, $\E=\Cat$, or the functor category $\Fun(I,\E)$ for an
essentially small $I$ and complete $E$. Since the categories $P$ and
$P^=$ are self-dual, Example~\ref{P.cone.exa} has a dual
counterpart: the image of a projector is also the limit, and taking
such images commutes with limits.

\subsection{Categories of elements.}

To compute the limit and colimit of a functor $X:I^o \to \Sets$, $I
\in \Cat$, define the {\em category of elements} $IX$ by
\begin{equation}\label{IX.eq}
IX = I /_\Y X,
\end{equation}
where $\Y$ is the Yoneda embedding \eqref{yo.eq}. Explicitly, $IX$
is the category of pairs $\langle i,x \rangle$, $i \in I$, $x \in
X(i)$, with morphisms $\langle i,x \rangle \to \langle i',x'
\rangle$ given by maps $f:i \to i'$ such that $X(f)(x') = x$. We
then have the forgetful functor $\pi:IX \to I$, $\langle i,x \rangle
\mapsto i$, and $\lim_{I^o}X = \Sec(I,IX)$, where the sections
category is actually discrete. In particular, for any functor
$i_\idot:J \to I$ from some essentially small $J$, functors
$x_\idot:J \to IX$ equipped with an isomorphism $\pi \circ x_\idot
\cong i_\idot$ are in one-to-one correspondence with sections $s \in
\Sec(J,i_\idot^*IX) \cong \lim_{J^o}i_\idot^{o*}X$, so that if the
colimit $i = \colim_Ji_\idot$ exists, and $X$ preserves the dual
limit $i = \lim_{J^o}i_\idot^o$, then for any $x_\idot:J \to IX$
such that $\pi \circ x_\idot \cong i_\idot$, $x = \colim_Jx_\idot
\in X(i) \subset IX$ exists and is preserved by $\pi$. Among other
things, this implies that if $I$ is Karoubi-closed, then so is $IX$.

For colimits, we have $\colim_{I^o}X = \pi_0(IX)$; the canonical
projection
\begin{equation}\label{colim.eq}
  IX \to \pi_0(IX) \cong \colim_{I^o}X
\end{equation}
decomposes $IX$ into a coproduct of its connected components
$(IX)_x$ numbered by elements $x \in \colim_{I^o}X$, and we have
$(IX)_x \cong IX_x$, for a decomposition $X \cong \copr_xX_x$ of $X$
into a disjoint union of functors $X_x:I^o \to \Sets$ with
$\colim_{I^o}X_x=\ppt$.

Alternatively, for a covariant functor $X:I \to \Sets$, we can
consider the {\em dual category of elements} $I^\perp X = (I^o /_Y
X)^o$. Explicitly, its objects are again pairs $\langle i,x
\rangle$, $i \in I$, $x \in X(i)$, and morphisms $\langle i,x
\rangle \to \langle i',x'\rangle$ are maps $f:i \to i'$ such that
$X(f)(x) = x'$. We also have the forgetful functor $\pi:I^\perp X
\to I$, and then $\lim_IX \cong \Sec(I,I^\perp X)$, $\colim_IX \cong
\pi_0(I^\perp X)$.

\begin{exa}\label{Y.exa}
For any $i \in I$, with the Yoneda image $\Y(i) \in I^o\Sets$, we
have $I\Y(i) \cong I / i$; this category has a terminal object, so
$\pi_0(I\Y(i)) \cong \ppt$, and we also have $\colim_{I^o}\Y(i)
\cong \ppt$.  Conversly, a functor $X \in I^o\Sets$ such that $IX$
has a terminal object $o$ lies in the essential image of
\eqref{yo.eq} (namely, we have $X \cong \Y(\pi(o))$.
\end{exa}

\begin{exa}\label{tw.exa}
Consider the $\Hom$-pairing $Y:I^o \times I \to \Sets$ for a small
category $I$. Then the category of elements $(I \times I^o)Y$ is the
{\em twisted arrows category} $\Tw(I)$; it objects are arrows $f:i
\to i'$ in $I$, and morphisms from $f_0:i_0 \to i_0'$ to $f_1:i_1
\to i_1'$ are given by commutative diagrams
\begin{equation}\label{tw.dia}
\begin{CD}
i_0 @>{f_0}>> i_0'\\
@VVV @AAA\\
i_1 @>{f_1}>> i_1'.
\end{CD}
\end{equation}
\end{exa}

\begin{exa}\label{IX.exa}
For any $i \in IX$, the comma-fiber $i \setminus IX$ of the
forgetful functor $\pi:IX \to I$ decomposes as $i \setminus IX \cong
\copr_{x \in X(i)}\langle i,x \rangle \setminus IX$. Then by
\eqref{kan.eq}, for any $Y:(IX)^o \to \Sets$, the left Kan extension
$\pi^o_!Y$ is given by $\pi^o_!Y(i) \cong \copr_{x \in
  X(i)}Y(\langle i,x \rangle)$. In particular, $\pi^o_!\ppt \cong
X$, and $\pi^o_!$ is conservative, so that the map $\pi^o_!Y \to X$
induced by $Y \to \ppt$ is an isomorphism if and only if $Y \to
\ppt$ is an isomorphism. We also have
\begin{equation}\label{yo.X.eq}
X \cong \colim_{\langle i,x \rangle \in IX}\Y(i);
\end{equation}
if $X = \ppt$, this is the isomorphism $\colim_IY \cong \ppt$, and
the general case is obtained by applying $\pi^o_!$ and
\eqref{yo.sq}.
\end{exa}

More generally, for any functor $\C:I^o \to \Cat$, we can define
a category $\ra{I\C}$ resp.\ $\la{I\C}$ as the category of pairs
$\langle i,c \rangle$, $i \in I$, $c \in \C(i)$, with morphisms
$\langle i,c \rangle \to \langle i',c' \rangle$ given by pairs
$\langle f,g \rangle$ of a morphism $f:i \to i'$ and a morphism $g:c
\to \C(f)(c')$ resp.\ $g:\C(f)(c') \to c$, and then define the {\em lax}
and {\em co-lax limit} of the functor $\C$ by
\begin{equation}\label{lax.eq}
\ra{\lim}_{I^o}\C = \Sec(I,\ra{I\C}), \qquad \la{\lim}_{I^o}\C =
\Sec(I^o,(\la{I\C})^o).
\end{equation}
If $\C \cong X$ takes values in $\Sets \subset \Cat$, then $\ra{I\C}
\cong \la{I\C} \cong IX$, and the lax and co-lax limits coincide
with the usual limit $\lim X$. For any $\C:I^o \to \Cat$, we have
the relative $\Hom$-pairing $\la{I\C} \times_I \ra{I\C} \to \Sets$,
and this induces a relative version
\begin{equation}\label{yo.lax}
\Y_I:\ra{\lim}_{I^o}\C \to \Fun(\la{I\C},\Sets)
\end{equation}
of the Yoneda embedding \eqref{yo.eq} sending a section $s:I \to
\ra{I\C}$ to the functor $\langle i,c \rangle \mapsto \Hom(c,s(i))$,
where $\langle i,s(i) \rangle \in \ra{I\C}$ is the value of $s$ at
$i$. In terms of left Kan extensions, we have
\begin{equation}\label{yo.lax.kan}
\Y_I(s) = \Tw(s)_!\ppt,
\end{equation}
where $\Tw(s):\Tw(I) \to \la{I\C}$ sends an arrow $f:i \to i'$ in
the twisted arrow category $\Tw(I)$ of Example~\ref{tw.exa} to
$\langle i,\C(f)(i') \rangle$.

\begin{exa}\label{1.exa}
Let $[1]$ be the ``single arrow category'' with two objects $0$,
$1$, and a single non-identity map $0 \to 1$. Then a functor
$\C_\idot:[1] \to \Cat$ is given by two small categories $\C_0$,
$\C_1$ equipped a functor $\gamma:\C_0 \to \C_1$, and we have
$\la{\lim}_{[1]}\C_\idot \cong \C_0 /_\gamma \C_1$ and
$\ra{\lim}_{[1]}\C_\idot \cong \C_1 \setminus_\gamma \C_0$.
\end{exa}

A functor $\gamma:I \to I'$ between essentially small categories is
a {\em localization} if for any category $\E$, $\gamma^*:\Fun(I',\E)
\to \Fun(I,\E)$ is fully faithful. If $I'=\ppt$, then $\gamma:I \to
\ppt$ is a localization iff $I$ is connected. For any localization
$\gamma:I \to I'$ and functor $E:I' \to \E$, the left Kan extension
$\gamma_!\gamma^*E$ exists and is given by $E$, so that we have a
natural isomorphism
\begin{equation}\label{cof.eq}
  \colim_I\gamma^*E \cong \colim_{I'}E
\end{equation}
whose source exists iff so does its target. More generally, a
functor $\gamma:I \to I'$ between essentially small categories is
{\em cofinal} if $i \setminus I$ is connected for any $i \in I'$, or
equivalently, by \eqref{kan.eq}, if $\gamma^o_!\ppt \cong \ppt$. In
this case, for any target category $\E$ and functor $E:I' \to \E$, we
still have the isomorphism \eqref{cof.eq} (if $I'=\ppt$, then
$\gamma$ is a localization, and the general case reduces to this by
\eqref{kan.eq}). For example, any functor $\gamma$ that admits a
left-adjoint is cofinal (in this case, all the comma-fibers $i
\setminus I$ have initial objects). Dually, $\gamma$ is {\em final}
iff $\gamma^o$ is cofinal, and for any final functor $\gamma$, we
have \eqref{cof.eq} with $\colim$ replaced by $\lim$. The opposite
of a localization is trivially a localization, so a localization is
both final and cofinal.

\begin{exa}\label{cof.exa}
If we are given a functor $\gamma:I' \to I$ between small
categories, and functors $X:I^o \to \Sets$, $X':{I'}^o \to \Sets$,
then a map $a:\gamma^o_!X' \to X$ is adjoint to a map $a^\dg:X' \to
\gamma^{o*}X$, thus defines a functor $\alpha:I'X' \to IX$.
Conversely, $a$ is recovered from $\alpha$ by applying the
conservative functor $\pi^o_!$ of Example~\ref{IX.exa} to the
tautological map $\alpha^o_!\ppt \to \ppt$, so that $\alpha$ is
cofinal iff $a$ is an isomorphism.
\end{exa}

\subsection{Partially ordered sets.}

A partially ordered set $J$ defines a small category in the usual
way --- objects are elements $j \in J$, and there is exactly one map
$j \to j'$ iff $j \leq j'$. We denote the category of partially
ordered sets by $\Pos$, and we have a full embedding $\Pos \subset
\Cat$, so it makes sense to say that a small category $\C$ ``is a
partially ordered set'' (explicitly, this means that there is at
most one morphism between any two objects $c,c' \in \C$, and the
only invertible morphisms in $\C$ are identity maps). In particular,
a discrete category is a partially ordered set, so that $\Pos
\supset \Sets$. For any $J \in \Pos \subset \Cat$, we have
$J^>,J^<,J^o \in \Pos \subset \Cat$. A {\em full subset} $J' \subset
J$ is a subset equipped with the induced order (so that the
embedding $J' \to J$ is full as a functor). The {\em fiber} $J_j$ of
a morphism $f:J \to J'$ in $\Pos$ is the full subset $f^{-1}(j) \in
J$, and the {\em left} and {\em right comma-fibers} $J/j$ resp.\ $j
\setminus J$ are the full subsets of elements $j' \in J$ such that
$f(j') \leq j$ resp.\ $f(j') \geq j$ (up to an equivalence, this
agrees with our general categorical usage). We denote $J /' j = J/j
\ssetminus \{j\}$, and we note that $J/j \cong (J /' j)^>$. A
morphism $f:J \to J'$ in $\Pos$ is {\em conservative} if it is
conservative as a functor, or equivalently, if all the fibers $J_j$,
$j \in J'$ are discrete.

\begin{exa}
For any $J \in \Pos$ and functor $X:J^o \to \Sets$, the category of
elements $JX$ is a partially ordered set, and the forgetful functor
$JX \to J$ is conservative.
\end{exa}

For any integer $n \geq 0$, we let $[n]$ be the set $\{0,\dots,n\}$
with the usual total order, so that a chain \eqref{c.idot} in some
$J \in \Pos$ is the same thing as a map $[n] \to J$. We say that
$\dim J \leq n$ if all the chains of length $> n$ are degenerate, or
equivalently, if $J$ admits a conservative map $J \to [n]$ (for
example, one can take the height map $\hht$ sending $j \in J$ to the
maximal length of a non-degenerate chain in $J / j$). In particular,
$[1]$ is the single arrow category of Example~\ref{1.exa}, and $\dim
J \leq 1$ iff there exists a conservative map $J \to [1]$.
        
A full subset $J_0 \subset J$ is {\em left-closed} if $J_0 / j
\subset J_0$ for any $j \in J_0$, and {\em right-closed} if $J_0^o
\subset J^o$ is left-closed. The complement $J_1 = J \ssetminus J_0
\subset J_0$ of a left-closed subset $J_0 \subset J$ is
right-closed, and vice versa. For any $J \in \Pos$, $J \subset J^>$
is left-closed, and $J \subset J^<$ is right-closed. For any set
$S$, we denote by $P(S) \in \Pos$ the set of all subsets $S' \subset
S$, ordered by inclusion, and for any $J \in \Pos$, we let $L(J)
\subset P(J)$ be the full subset of left-ordered full subsets. The
embedding $L(J) \to P(J)$ admits a left-adjoint
\begin{equation}\label{lp.eq}
\Lambda:P(J) \to L(J), \qquad S \mapsto \bigcup_{s \in S}J/s.
\end{equation}
For any map $f:J \to J'$ in $\Pos$, we have a map $f_\dg:P(J) \to
P(J')$ sending $J_0 \subset J$ to $f(J_0) \subset J'$, and the
right-adjoint map $f^\dg:P(J') \to P(J)$, $J_0 \mapsto
f^{-1}(J_0)$. The adjoint map $f^\dg$ sends $L(J') \subset P(J')$
into $L(J) \subset P(J)$. The $\Hom$-pairing $J^o \times J \to
\Sets$ factors through the essential image of the fully faithful
embedding $[1] \to \Sets$, $0 \mapsto \varnothing$, $1 \mapsto
\ppt$, we have $J^o[1] \cong L(J)$, and the Yoneda embedding
\eqref{yo.eq} reduces to a map $\Y:J \to L(J)$ sending $j \in J$ to
$J/j \subset J$. This leads to the following universal gluing
construction for left-closed embeddings.

\begin{exa}\label{glue.exa}
For any $J_0,J_1 \in \Pos$ and a map $\lambda:J_1 \to L(J_0)$,
define the {\em gluing} $J= J_0 \copr_\lambda J_1$ as the coproduct
$J_0 \copr J_1$, with the order $j \leq j'$ iff either $j,j' \in
J_l$, $l=0,1$, $j \leq j'$, or $j' \in J_1$, $j \in \lambda(j')
\subset J_0$. Then $J_0 \subset J$ is left-closed, $J_1 = J
\ssetminus J_0 \subset J$ is right-closed, and if we let $\eps_l:J_l
\to J$ be the embeddings, then $\lambda = \eps_0^\dg \circ \Y \circ
\eps_1$. Conversely, for any left-closed full embedding $\eps_0:J_0
\to J$ in $\Pos$, with the complementary right-closed full embedding
$\eps_1:J_1 \to J$, we have $J \cong J_0 \copr_\lambda J_1$ for
$\lambda = \eps_0^\dg \circ \Y \circ \eps_1$.
\end{exa}

While $\Pos$ is cocomplete, colimits in $\Pos$ often behave badly
--- in particular, they need not be preserved by the forgetful
functor $\Pos \to \Sets$. The gluing construction of
Example~\ref{glue.exa} shows that this does not happen for pushouts
of left-closed embeddings. Namely, for any left-closed full
embedding $J_0 \to J \cong J_0 \copr_\lambda J_1$, and any map
$f:J_0 \to J'$ to some $J' \in \Pos$, we have a cocartesian square
\begin{equation}\label{st.sq}
\begin{CD}
J_0 @>{f}>> J'\\
@VVV @VVV\\
J @>>> J' \copr_{\lambda'} J_1,
\end{CD}
\end{equation}
where $\lambda' = \Lambda \circ f_\dg \circ \lambda$, and $\Lambda$
is the map \eqref{lp.eq}. In particular, for any $J \in \Pos$ of
dimension $\dim J = n$, we can choose a conservative map $J \to
[n]$, with discrete fiber $J_n$ and a left comma-fiber
$J_{<n}=J/(n-1)$ of dimension $\dim J_{<n} \leq n-1$; then $J \cong
J_{<n} \copr_\lambda J_n$ for a map $\lambda:J_n \to L(J_{<n})$
sending $j \in J_n$ to $J_{<n} / j = J /' j$, and \eqref{st.sq}
provides a cocartesian square
\begin{equation}\label{dim1.sq}
\begin{CD}
\coprod_{j \in J_n}\lambda(j) @>>> \coprod_{j \in J_1}\lambda(j)^>\\
@VVV @VVV\\
J_{<n} @>>> J.
\end{CD}
\end{equation}
where the map on the left is the coproduct of the embeddings
$\lambda(j) \subset J_{<n}$. This allows to build up
finite-dimensional partially ordered sets by induction on dimension.

For any regular cardinal $\kappa$ and $J \in \Pos$, $J$ is
$\kappa$-bounded as a set iff it is $\kappa$-bounded as a small
category, and we let $\Pos_\kappa \subset \Pos$ be the full
subcategory spanned by $\kappa$-bounded $J \in \Pos$. We say that $J
\in \Pos$ is {\em left-$\kappa$-bounded} if $|J/j| < \kappa$ for all
$j \in J$, and {\em left-finite} if it is left-$\kappa$-bounded for
the countable cardinal $\kappa$ (that is, $J/j$ is finite for any
$j$). We let $P_\kappa(J) \subset P(J)$, $L_\kappa(J) \subset L(J)$
be the full subsets spanned by $\kappa$-bounded subsets $J' \subset
J$, and we note that $J$ is left-$\kappa$-bounded iff \eqref{lp.eq}
sends $P_\kappa(J)$ into $L_\kappa(J)$.

\subsection{Nerves and replacements.}

As usual, we let $\Delta \subset \Pos \subset \Cat$ be the full
subcategory spanned by $[n]$, $n \geq 0$. Functors $\Delta^o \to
\Sets$ are {\em simplicial sets}, and the {\em nerve functor}
$N:\Cat \to \Delta^o\Sets$ is given by $N = \phi^o_!\Y$, where
$\phi:\Delta \to \Cat$ is the embedding, and $\Y$ is the Yoneda
embedding \eqref{yo.eq} for $I = \Delta$. If one computes $\phi^o_!$
by \eqref{kan.eq}, one recovers the usual definition of the nerve:
$N(\C)([n])$ is the set of chains \eqref{c.idot} in a small category
$\C$. Note that sending a chain to its last element $c_n$ provides a
well-defined functor
\begin{equation}\label{xi.eq}
  \Delta N(\C) \to \C,
\end{equation}
where $\Delta N(\C)$ is the category of elements of the nerve
$N(\C):\Delta^o \to \Sets$. It is also useful to note that the nerve
functor $N$ is fully faithful. In particular, it is easy to see that
it sends a cocartesian square \eqref{dim1.sq} to a cocartesian
square in $\Delta^o\Sets$, and this implies that \eqref{dim1.sq} is
cocartesian in $\Cat \supset \Pos$.

In principle, one can recover a small category $\C$ from its nerve
$N(\C)$; however, for many applications, the full nerve is not
needed. Namely, let $\V = \{0,1\}^<$ be the partially ordered set
with three elements $o,0,1$ and order $o \leq 0,1$, and let
$\mu:\V^o \to \Delta$ be the functor sending $0,1$ to $[0]$, $o$ to
$[1]$, and the maps $0,1 \to o$ to the embeddings $[0] \to [1]$ onto
$0,1 \in [1]$. Then for any simplicial set $X:\Delta^o \to \Sets$,
the category of elements $\V^o\mu^*X$ is a partially ordered set of
dimension $\dim \V^o\mu^*X \leq 1$ equipped with a functor
$\mu:\V^o\mu^*X \to \Delta X$, and for any small category $\C$, the
corresponding partially ordered set $\V(\C) = \V^o\mu^*N(\C)$ comes
equipped with a functor
\begin{equation}\label{V.xi.eq}
q=\xi \circ \mu:\V(\C) \to \C
\end{equation}
induced by the functor \eqref{xi.eq}. For any small $\C$, $\V(\C)$
is left-finite, we have $\dim \V(\C) \leq 1$, and if $\C$ is
$\kappa$-bounded for some regular $\kappa$, then so is
$\V(\C)$. Note that we have a canonically isomorphism $\V(\C) \cong
\V(\C^o)$ (if $\C=\ppt$, this is the isomorphism $\V^o \to \V^o$
interchanging $0$ and $1$).

\begin{lemma}\label{V.le}
For any small category $\C$, the functor \eqref{V.xi.eq} is a
localization.
\end{lemma}

\proof{} By definition, for any functors $E,E':\C \to \E$ to some
category $\E$, a map $\alpha:E \to E'$ is a collection of maps
$\alpha(c):E(c) \to E'(c)$ such that for any map $f:c \to c'$ in
$\C$, we have $E'(f) \circ \alpha(c) \cong \alpha(c') \circ
E(f)$. Equivalently, we might take a collection of pairs of maps
$\alpha(c)_0,\alpha_1(c):E(c) \to E(c')$ such that $E'(f) \circ
\alpha(c)_0 \cong \alpha(c')_1 \circ E(f)$ -- indeed, this condition
for $\id:c \to c$ insures that $\alpha(c)_0 = \alpha(c)_1$, so the
data are exactly the same. But then these data can be packaged into
a single limit over $\V(C)$, and this provides a canonical
isomorphism
\begin{equation}\label{V.eq}
\Hom(E,E') \cong \lim_{v\in\V(\C)}\Hom(E(q_\perp(v)),E'(q(v))),
\end{equation}
where $q$ is the functor \eqref{V.xi.eq}, and $q_\perp:\V(\C) \to
\C^o$ is obtained by identifying $\V(\C)\cong\V(\C^o)$ and taking
\eqref{V.xi.eq} for $\C^o$. it remains to observe that
$\Hom(q^*E,q^*E')$ is given by exactly the same expression.
\endproof

It is useful to note that one can further refine Lemma~\ref{V.le} by
describing the essential image of the full embedding
$q^*:\Fun(\C,\E) \to \Fun(\V(\C),\E)$ --- namely, it consists of
functors $\V(\C) \to \E$ that invert all maps in $\V(\C)$ lying over
the arrow $0 \to o$ in $\V^o$. More formally, let $e:[1] \to \V^o$
be the embedding onto this arrow, $0 \mapsto 0$, $1 \mapsto o$, and
let $\V(\C)_0 \subset \V(\C)$ be the partially ordered set $\V(\C)_0
= [1]e^*\mu^*N(\C) \subset \V(\C)$. Then if we let $\C_0 =
N(\C)([0])$ be the set of objects of $\C$, considered as a discrete
category, we have a projection $p:\V(\C)_0 \to \C_0$ right-adjoint
to the embedding $\C_0 \to \V(\C)_0$, and we then have a cartesian
square
\begin{equation}\label{V.sq}
\begin{CD}
\Fun(\C,\E) @>{q^*}>> \Fun(\V(\C),\E)\\
@V{s}VV @VV{e^*}V\\
\Fun(\C_0,\E) @>{p^*}>> \Fun(\V(\C)_0,\E),
\end{CD}
\end{equation}
where $t:\C_0 \to \C$ is the tautological embedding.

\begin{corr}\label{V.corr}
For any regular cardinal $\kappa$, a category $\E$ is
$\kappa$-complete resp.\ $\kappa$-cocomplete iff for any $J \in
\Pos_\kappa$ of dimension $\dim J \leq 1$, and any functor $E:J \to
\E$, $\lim_JE$ resp.\ $\colim_JE$ exists, and for any two
$\kappa$-complete resp.\ $\kappa$-cocomplete categories $\E$, $\E'$,
a functor $F:\E \to \E'$ preserves $\kappa$-bounded limits
resp.\ colimits iff it preserves the limits $\lim_J$ resp.\ the
colimits $\colim_J$ for all $J \in \Pos_\kappa$, $\dim J \leq 1$.
\end{corr}

\proof{} Clear. \endproof

\section{Filtered colimits.}

Let us now introduce the first of our two main definitions. Fix a
regular cardinal $\kappa$.

\begin{defn}\label{filt.def}
A functor $\gamma:\C \to \C'$ between essentially small categories
is {\em $\kappa$-filtered} if for any commutative square
\begin{equation}\label{filt.sq}
\begin{CD}
I @>{a}>> \C\\
@V{\eps}VV @VV{\gamma}V\\
I^> @>{b}>> \C'
\end{CD}
\end{equation}
with $\kappa$-bounded small $I$, there exists a functor $q:I^> \to
\C$ such that $q \circ \eps \cong a$ and $\gamma \circ q \cong
b$. An essentially small category $\C$ is {\em $\kappa$-filtered} if
so is the tautological projection $\C \to \ppt$ (or equivalently, if
any functor $I \to \C$ from a $\kappa$-bounded small category $I$
admits a cone $I^> \to \C$). A category $\C$ is {\em
  $\kappa$-f-cocomplete} if any functor $E:I \to \C$ from a
$\kappa$-filtered $I$ admits a colimit $\colim_IE$. An object $c
\in \C$ in a $\kappa$-f-cocomplete category $\C$ is {\em
  $\kappa$-compact} if $\Hom(c,-):\C \to \Sets$ preserves
$\kappa$-filtered colimits.
\end{defn}

\begin{remark}
The term ``filtered'' is actually a mistranslation of the original
French ``filtrant'' (the correct translation is
``filtering''). Unfortunately, it is too late to fix this.
\end{remark}

Compositions and pullbacks of $\kappa$-filtered functors are
trivially $\kappa$-filtered, and for any collection of
$\kappa$-filtered functors $\gamma_s:\C_s \to \C$ indexed by some
set $S$, their product
\begin{equation}\label{prod.rel}
  \prod_s\gamma_s:\C \times_{\C^S} \prod_{s \in S}\C_s \to \C
\end{equation}
is also $\kappa$-filtered (in particular, a product of
$\kappa$-filtered categories is $\kappa$-filtered). A
$\kappa$-filtered category is obviously connected. We denote by
$\Filt_\kappa \subset \Cat$ the full subcategory spanned by
$\kappa$-filtered small categories. For any regular cardinals
$\kappa' > \kappa$, a $\kappa'$-filtered category is
$\kappa$-filtered, and consequently, a $\kappa$-f-cocomplete
category $\C$ is $\kappa'$-f-cocomplete, and a $\kappa$-compact
object in such a category is $\kappa'$-compact. If $\kappa$ is the
countable cardinal, then $\kappa$-filtered, $\kappa$-f-cocomplete,
$\kappa$-compact are shortened to ``filtered'', ``f-cocomplete'',
``compact''. Any complete category is trivially $\kappa$-f-complete
for any regular $\kappa$ -- in particular, this includes $\Sets$ and
the functor category $I^o\Sets$ for any essentially small $I$. For
any $i \in I$, its Yoneda image $\Y(i) \in I^o\Sets$ is
$\kappa$-compact (indeed, $\Hom(\Y(i),-)$ is simply evaluation at
$i$).

\begin{lemma}\label{cof.le}
Assume given a $\kappa$-filtered category $\C$, and a fully faithful
embedding $\gamma:\C' \subset \C$ such that
$c \setminus \C'$ is non-empty for any $c \in \C$. Then $\C'$ is
$\kappa$-filtered, and $\gamma$ is cofinal.
\end{lemma}

\proof{} For any functor $E:I \to \C' \subset \C$ with a cone
$E_>:I^> \to \C$, with some vertex $c$, we have a map $f:c \to c'$
with $c' \in \C' \subset \C$, and then $f_!E_>:I^> \to \C'$ is a
cone for $E$. Therefore $\C'$ is $\kappa$-filtered. Moreover, for
any $c \in \C$, a functor $E:I \to c \setminus \C'$ is the same
thing as a functor $E^<:I^< \to \C$ sending $o \in I^<$ to $c$ and
$I \subset I^<$ into $\C' \subset \C$, and the same construction
applied to $E^<$ then provides a cone $E_>:I^> \to c \setminus \C'$
for $E$. Therefore $c \setminus \C'$ is also $\kappa$-filtered,
hence connected.
\endproof

\begin{lemma}\label{pos.filt.le}
Assume given an essentially small category $\C$ such that for any
left-finite $\kappa$-bounded partially ordered set $J \in
\Pos_\kappa$ of chain dimension $\dim J \leq 1$, any functor $J \to
\C$ admits a cone $J^> \to \C$. Then $\C$ is $\kappa$-filtered.
\end{lemma}

\proof{} By \eqref{cone.eq}, any cofinal functor $\gamma:I' \to I$
between essentially small categories induces an equivalence
$\Cone(\gamma^*E) \cong \Cone(E)$ for any $E:I \to \C$. It remains
to invoke Lemma~\ref{V.le}.
\endproof

\begin{exa}\label{pos.filt.exa}
Lemma~\ref{pos.filt.le} immediately shows that a partially ordered
set $J \in \Pos$ is $\kappa$-filtered iff for any $\kappa$-bounded
subset $J' \subset J$, $|J'| < \kappa$, there exist an upper bound
$j \in J$ (that is, $j \in J$ such that $J' \subset J/j \subset J$).
\end{exa}

\begin{exa}\label{LP.exa}
In particular, for any set $S$, the partially ordered set
$P_\kappa(S)$ of $\kappa$-bounded subsets $S' \subset S$ is
$\kappa$-filtered, and for a partially ordered set $J$, so is the set
$L_\kappa(J)$ of left-closed $\kappa$-bounded subsets $J'\subset
J$. In both cases, an upper bound is given by the union. Among other
things, this implies that a $\kappa$-compact set $S \in \Sets$ is
$\kappa$-bounded --- indeed, $S \cong \colim_{P_\kappa(S)}S'$, so
$\id:S \to S$ must factor through the inclusion $S' \to S$ for some
$S' \in P_\kappa(S)$.
\end{exa}

A small category $I$ such that the identity functor $\id:I \to I$
admits a cone $\id_>:I^> \to I$ is trivially $\kappa$-filtered for
any regular cardinal $\kappa$. In particular, this applies to any
category that has a terminal object -- take the universal cone --
and also to the category $P$ of Example~\ref{P.cone.exa}, so that
any $\kappa$-f-cocomplete category $\C$ is Karoubi-closed. These two
examples are essentially universal, as seen from the following
result.

\begin{lemma}\label{P.le}
Assume given an essentially small category $I$. Then $\id:I \to I$
admits a cone iff the Karoubi closure $P(I)$ of the category $I$ has
a terminal object.
\end{lemma}

\proof{} By definition, for any cone $\id_>$ for $\id:I \to I$, with
some vertex $o$, we have a map $p(i):i \to o$ for any $i \in I$ such
that $p(i') \circ f = p(i)$ for any map $f:i \to i'$. In particular,
$p=p(o):o \to o$ is an idempotent, $p^2=p$, and then the
corresponding object $\langle o,p \rangle$ in $P(I)$ is
terminal. Conversely, for any terminal object $\langle o,p \rangle$
in $P(I)$, we can let $p(i):i \to o$ be the composition of the
unique map $\langle i,\id \rangle \to \langle o,p \rangle$ and the
embedding $a:\langle o,p \rangle \to \langle o,\id \rangle$, and
this defines a cone $\id_>$.
\endproof

\begin{lemma}\label{con.le}
For any $\kappa$-filtered category $I$ and functor $X:I \to \Sets$,
the dual category of elements $I^\perp X$ is $\kappa$-filtered iff
$\colim_IX=\ppt$.
\end{lemma}

\proof{} Since $\colim_IX \cong \pi_0(I^\perp X)$, the ``only if''
part is clear. For the ``if'' part, assume given a functor $\gamma:J
\to I^\perp X$, $j \mapsto \langle i(j),x(j) \rangle$ from a
$\kappa$-bounded $J$. Since $I$ is $\kappa$-filtered, its
composition with the projection $I^\perp X \to I$ has a cone, with
some vertex $i \in I$ and maps $f(j):i(j) \to i$, $j \in J$, and the
functor $\gamma':J \to X(i) \subset I^\perp x$, $j \mapsto \langle
i,X(f(j))(x(j)) \rangle$ comes equipped with a map $\gamma \to
\gamma'$, so it suffices to find a cone for $\gamma'$. Since $X(i)$
is discrete, $\gamma'$ factors through the set $S = \pi_0(J)$, and
it then suffices to find a map $g:i \to i'$ such that $X(g) \circ
\gamma':S \to X(i')$ sends any two elements $s,s' \in S$ to the same
element in $X(i')$. For any given two elements $s,s'$, such an $i' =
i'(s,s')$ exists since $\colim_IX \cong \ppt$, and then since $S
\times S$ is $\kappa$-bounded, one can find a single $i'$ for all
pairs by considering the corresponding functor $(S \times S)^< \to
I$, $o \to i$, $s \times s' \to i'(s,s')$, and then taking its cone.
\endproof

\begin{corr}\label{con.corr}
For any cofinal functor $\gamma:I \to I'$ between $\kappa$-filtered
categories, and any $i \in I'$, the right comma-fiber $i \setminus
I$ is $\kappa$-filtered.
\end{corr}

\proof{} We have $i \setminus I \cong I^\perp X$, where the functor
$X:I \to \Sets$ sends $i' \in I$ to the set $\Hom(i,\gamma(i'))$.
\endproof

\begin{lemma}\label{cone.le}
For any $\kappa$-filtered category $\C$ and $\kappa$-bounded small
category $J$, the functor category $\Fun(J,\C)$ is
$\kappa$-filtered, and so is the cone category $\Cone(E)$ for any
functor $E:J \to \C$, while the embedding $\gamma^*:\C \to
\Fun(J,\C)$ induced by the tautological projection $J \to \ppt$ is
cofinal.
\end{lemma}

\proof{} For the first claim, note that a functor $I \to \Fun(J,\C)$
is the same thing as a functor $I \times J \to \C$, and to construct
a cone for it, it suffices to take a cone $(I \times J)^> \to \C$
and compose it with the projection $I^> \times J \to (I \times J)^>$
induced by \eqref{bicone.eq}. For the second claim, a functor $F:I \to
\Cone(E)$ is the same thing as a functor $F':I \times J^> \to \C$ whose
restriction to $I \times J$ is identified with $(\id \times
\gamma)^*E$; to obtain a cone for $F$, consider the product $I
* J$ of \eqref{star.eq}, extend $F'$ to a functor $F'':I * J \to \C$
by letting $F'' = E$ on $\{o\} \times J \subset I * J$, and take a
cone of $F''$. For the third claim, note that by \eqref{cone.eq},
all the comma-fibers $E \setminus_{\gamma^*} \C$ are
$\kappa$-filtered, thus connected.
\endproof

\begin{lemma}\label{dim1.le}
For any left-$\kappa$-bounded partially ordered set $J$ of some
finite dimension $\dim J \leq n$, the co-lax limit $\la{\lim}_JX$ of
any functor $X:J \to \Filt_\kappa$ is $\kappa$-filtered.
\end{lemma}

\proof{} Let $\C = \la{\lim}_JX$. If $\dim J=0$, so that $\C \cong
\prod_{j \in J}X(j)$, the claim is clear. In general, by induction
on $n$, we may choose a conservative map $J \to [n]$, and assume the
claim proved for $J_{<n} = J/(n-1) \subset J$. Then it suffice to
prove that the restriction functor $\la{\lim}_JX \to
\la{\lim}_{J_{<n}}X$ is $\kappa$-filtered. By \eqref{dim1.sq}, this
functor is a product \eqref{prod.rel} of pullbacks of the
restriction functors for the embeddings $J /' j \subset J/j$, $j \in
J_n$, so it suffices to assume that $J = {J'}^>$ for some
$\kappa$-bounded $J' \in \Pos_\kappa$, and prove that the functor
$\eps^*:\C \to \C' = \la{\lim}_{J'}X$ is $\kappa$-filtered. Indeed,
for any $j \in J'$, we have the functor $X(j) \to X(o)$, so that an
object $c \in \C'$ defines a functor $c_\dg:J' \to X(o)$, and then
the fiber $\C_c$ of the functor $\eps^*$ is the cone category
$\Cone(c_\dg)$. Then a diagram \eqref{filt.sq} defines a functor $J'
* I \to X(o)$, where $*$ is the product \eqref{star.eq}, and to find
the required $q$, it suffices to construct a cone for this
functor. Since $J'$ and $I$ are $\kappa$-bounded and $X(o)$ is
$\kappa$-filtered, this can be done.
\endproof

\begin{corr}\label{comma.corr}
For any cofinal functor $\gamma:\C \to I$ between $\kappa$-filtered
categories $\C$, $I$, the full embedding $\eta:\C \to \C /_\gamma I$
of \eqref{comma.facto} is cofinal.
\end{corr}

\proof{} By Example~\ref{1.exa} and Lemma~\ref{dim1.le}, $\C
/_\gamma I$ is $\kappa$-filtered, so by Lemma~\ref{cof.le}, it
suffices to find an object in $\langle c,i,\alpha \rangle \setminus_\eta
\C$ for any $\langle c,i,\alpha \rangle \in \C /_\gamma
I$. Indeed, by Corollary~\ref{con.corr}, the comma-fiber $i
\setminus_\gamma \C$ is filtered, so we have an object $c' \in \C$
and a map $f:i \to \gamma(c')$. We then have two objects $\langle
c,\id \rangle$, $\langle c',f \circ \alpha \rangle$ in $\gamma(c)
\setminus_\gamma \C$, and this comma-fiber is also filtered, so we
can choose $c'' \in \C$ and maps $g:c \to c''$, $g':c' \to c''$ such
that $\gamma(g) = \gamma(g') \circ f \circ \alpha$. Then $c''$ with
the map $\langle g,\gamma(g') \circ f\rangle:\langle c,i,\alpha
\rangle \to \langle c'',\gamma(c''),\id \rangle$ is in
$\langle c,i,\alpha \rangle \setminus_\eta \C$.
\endproof

\begin{prop}\label{filt.prop}
An essentially small category $I$ is $\kappa$-filtered iff
$\colim_I$ commutes with $\kappa$-bounded limits in $\Sets$.
\end{prop}

\proof{} For the ``if'' part, note that by \eqref{cone.eq}, for any
functor $E:J \to I$ from an essentially small $J$, the set
$\Cone(E,i)$ of cones $E_>:J^> \to I$ with some fixed vertex $i$ is
given by $\lim_{j \in J^o}\Hom(E(j),i)$. If $J$ is
$\kappa$-bounded, we have
$$
\colim_{i \in I}\Cone(E,i) \cong \lim_{j \in
  J}\colim_{i \in I}\Hom(E(j),i) \cong \lim_{j \in J}\ppt \cong
\ppt
$$
by assumption, so that $\Cone(E,i)$ is non-empty at least for one $i
\in I$.

For the ``only if'' part, note that by Corollary~\ref{V.corr}, it
suffices to check that $\colim_I$ preserves $\lim_J$ for a
$\kappa$-bounded partially ordered set $J$ of dimension $\dim J \leq
1$. Fix such a $J$, assume given a functor $X:I \times J \to \Sets$,
and let $Y = \lim_JX:I \to \Sets$, $Z = \colim_IX:J \to
\Sets$. Moreover, let $Y' = \lim_J(\ev \times \id)^*X:\Fun(J,I)
\times J \to \Sets$, where $\ev$ is the evaluation pairing of
\eqref{fun.facto}. Then the map $\colim_IY \to \lim_JZ$ factors as
\begin{equation}\label{cofun.facto}
\begin{CD}
\colim_IY @>{\alpha}>> \colim_{\Fun(J,I)}Y' @>{\beta}>> \lim_JZ,
\end{CD}
\end{equation}
where $\alpha$ is induced by the tautological embedding $I \to
\Fun(J,I)$, thus invertible by Lemma~\ref{cone.le}. For $\beta$,
consider the categories of elements $\Z = (J^oZ)^o$, $\X = ((I
\times J)^oX)^o$, $\YY' = (\Fun(J,I)^oY')^o$, and note that $\YY'
\cong \Sec(J,\X)$, where $\X$ is considered as a category over $J$
via the canonical projection $\pi:\X \to \Z$ of \eqref{colim.eq} and
the projection $\Z \to J$. Then $\pi$ has $\kappa$-filtered fibers
by Lemma~\ref{con.le}, so that $\Sec(J,\pi):\Sec(J,\X) \to
\Sec(J,\Z) \cong \lim_JZ$ has $\kappa$-filtered fibers by
Lemma~\ref{dim1.le}, and $\beta=\pi_0(\Sec(J,\pi))$ is invertible.
\endproof

\begin{remark}
Our proof of Proposition~\ref{filt.prop} actually gives slightly
more. Namely, the fact that $\beta$ in \eqref{cofun.facto} is an
isomorphism only requires $J$ to be left-$\kappa$-bounded, not
necessarily $\kappa$-bounded. In particular, this holds for the
partially ordered set $\V(\C)$ for {\em any} small category $\C$
(and any regular $\kappa$). If one replaces $\Sets$ by an abelian
category $\A$, then Proposition~\ref{filt.prop} is Grothendieck's
axiom $AB5$ of \cite{toho}, while the fact that $\beta$ is an
isomorphism is the additional axiom $AB6$ (our proof easily implies
that both are satisfied by the category of abelian groups).
\end{remark}

\begin{corr}\label{filt.corr}
For any essentially small $I$, a functor $X:I^o \to \Sets$ is
$\kappa$-compact in $I^o\Sets$ iff there exists a $\kappa$-bounded
left-finite partially ordered set $J \in \Pos_\kappa$ of dimension
$\dim J \leq 1$, and a functor $\gamma:J \to I$ such that $X$ is a
retract of $\gamma^o_!\ppt$.
\end{corr}

\proof{} The ``if'' part immediately follows from
Proposition~\ref{filt.prop}: the functor $\Hom(\gamma^o_!\ppt,-)$ is
a limit over $J$, thus commutes with $\kappa$-filtered colimits. For
the ``only if'' part, consider the category $IX$ with the projection
$\pi:IX \to I$, so that $X \cong pi^o_!\ppt$, let $J_0 = \V(IX)$,
with the final functor $\xi:J \to IX$ of Lemma~\ref{V.le}, and let
$\gamma_0=\xi \circ \pi:J_0 \to I$, so that $X \cong \pi^o_!\ppt
\cong \pi^o_!\xi^o_!\ppt \cong \gamma_{0!}^o\ppt$. Then $J_0$ is
left-finite, $L_\kappa(J_0)$ is $\kappa$-filtered, and by
\eqref{kan.eq}, $\ppt:J_0^o \to \Sets$ coincides with $\colim_{J \in
  L_\kappa(J_0)} \eps(J)^o_!\ppt$, where $\eps(J):J \to J_0$ is the
embedding. Then since $X$ is compact, $\id:X \to X \cong \colim_{J
  \in L_\kappa(J_0)}\gamma_{0!}^o\eps(J)^o_!\ppt$ factors through
some $\gamma_{0!}^o\eps(J)^o_!\ppt$, and it suffices to take $\gamma
= \gamma_0 \circ \eps(J)$.
\endproof

\begin{corr}\label{sets.corr}
A set $S \in \Sets$ is $\kappa$-compact iff it is $\kappa$-bounded.
\end{corr}

\proof{} The ``only if'' part is Example~\ref{LP.exa}, and the
``if'' part is Proposition~\ref{filt.prop} for $J=S$.
\endproof

\begin{corr}\label{sets.ka.corr}
For any regular cardinals $\kappa' > \kappa$, the category
$\Sets_\kappa$ is $\kappa'$-f-cocomplete.
\end{corr}

\proof{} Since $\Sets$ is cocomplete, it suffices to show that for
any $X:I \to \Sets_\kappa$ with a $\kappa'$-filtered $I$,
$\colim_IX$ is $\kappa$-bounded. If not, there exists an injective
map $S \to \colim_IX$ with $|S|=\kappa < \kappa'$, and then by
Corollary~\ref{sets.corr}, it must factor through $X(i) \in
\Sets_\kappa$ for some $i \in I$. This is a contradiction.
\endproof

\begin{remark}
Note that for $\kappa'=\kappa$, Corollary~\ref{sets.ka.corr} is
certainly {\em not} true --- e.g. a filtered colimit of finite sets
can be as large as you want.
\end{remark}

\section{Directed posets.}\label{pos.sec}

A filtered partially ordered set $J$ is also called {\em directed};
by Example~\ref{pos.filt.exa} and an induction on cardinality, $J$
is directed iff for any elements $j_0,j_1 \in J$, there exists $j
\in J$ such that $j_0,j_1 \leq j$. Standard examples of directed
partially ordered sets are the set $\N$ of non-negative integeres,
with the standard order, and the set of finite subsets of a set $S$,
ordered by inclusion. More generally, we can also consider the set
of all finite left-closed subsets in a left-finite partially ordered
set $J$, again ordered by inclusion. A $\kappa$-filtered
generalization of the latter two examples is given by
Example~\ref{LP.exa}, and the first one is generalized as follows.

\begin{exa}\label{ord.exa}
Assume given a regular cardinal $\kappa$. Take a set $S_\kappa$ of
cardinality $\kappa$, and choose an order making it into a
well-ordered set $Q_\kappa$. Assume that the order is minimal in the
sense that $Q_\kappa$ is left-$\kappa$-bounded --- this can always
be achieved by taking the smallest $q \in Q_\kappa$ such that
$|Q_\kappa/q| = \kappa$, and replacing $Q_\kappa$ with
$Q_\kappa/'q$. Then $Q_\kappa$ is $\kappa$-filtered. Indeed, since
$\kappa$ is regular, the union $\cup_{s \in S'}Q_\kappa/s \subset
Q_\kappa$ for any $\kappa$-bounded subset $S' \subset Q_\kappa$ is
still $\kappa$-bounded, so it cannot cover the whole $Q_\kappa$.
\end{exa}

Since by Lemma~\ref{P.le}, a $\kappa$-bounded partially ordered set
is $\kappa$-filtered only if it has the largest element, the set $Q$
of Example~\ref{ord.exa} is the smallest non-trivial
$\kappa$-filtered partially ordered set. On the other hand, the sets
of Example~\ref{LP.exa} are rather big. To control cardinality of
$\kappa$-filtered partially ordered sets, the following notion turns
out to be very useful.

\begin{defn}\label{sharp.def}
For any regular cardinals $\kappa < \mu$, we have $\kappa \trl
\mu$ iff for any $\mu$-bounded set $S$, the partially ordered set
$P_\kappa(S)$ contains a cofinal $\mu$-bounded subset $U$.
\end{defn}

\begin{exa}\label{ka.ka.exa}
  For any regular cardinals $\kappa$, $\kappa'$, $\kappa''$ such
  that $\kappa \trl \kappa'$ and $\kappa' \trl \kappa''$, we have
  $\kappa \trl \kappa''$ (for any $\kappa''$-bounded set $S$, choose
  a cofinal $\kappa''$-bounded cofinal subset $J \subset
  P_{\kappa'}(S)$, and then choose a $\kappa'$-bounded cofinal
  subset in $P_\kappa(X)$ for any $X \in J$).
\end{exa}

\begin{exa}\label{suc.exa}
For any regular cardinal $\kappa$ with successor cardinal
$\kappa^+$, we have $\kappa \trl \kappa^+$. Indeed, for any $S$ with
$|S|=\kappa$, we can equip $S$ with a minimal order $Q_\kappa$ of
Example~\ref{ord.exa}, and then the embedding $Q_\kappa \to
P_\kappa(S)$, $q \mapsto Q_\kappa/q$ is cofinal.
\end{exa}

\begin{exa}\label{P.ka.exa}
For any regular cardinal $\kappa$, we have $\kappa \trl
P(\kappa)^+$. Indeed, for any set $S$ with $|S|=\kappa$ and any set
$S'$, a map $S' \to P(S)$ is the same thing as a subset in $S'
\times S$, and $|S' \times S| = |S|$ are soon as $S'$ is non-empty
and $\kappa$-bounded, so that $|P_\kappa(P(S))| = |P(S)|$. More
generally, for any $\kappa' \leq \kappa$, we have
$|P_{\kappa'}(P(S))| \leq |P_\kappa(P(S))| < P(\kappa)^+$, so that
$\kappa' \trl P(\kappa)^+$.
\end{exa}

\begin{exa}\label{sharp.non.exa}
For any uncountable regular cardinal $\kappa$, consider the
non-regular cardinal $\kappa_\infty$ of Example~\ref{reg.exa}, with
successor cardinal $\kappa_\infty^+$. Then it is {\em not} true that
$\kappa \trl \kappa_\infty^+$. Indeed, take a set $S$ with
$|S|=\kappa_\infty$, and for any subset $J \in P_\kappa(S)$, let
$U(J) = \cup_{X \in J} \subset S$ be the union of all subsets $X \in
J$. Then $|U(J)| \leq \max(\kappa,|J|)$. If $J \subset P_\kappa(S)$
is cofinal, we have $U(J)=S$, so that if $J$ is
$\kappa_\infty^+$-bounded, the only option is
$|J|=\kappa_\infty$. Represent $J$ as a union $J = \cup_nJ_n$ of
subsets $J_0 \subset \dots \subset J_n \subset \dots$ of
cardinalities $|J_n|=\kappa_n$. Then $S = \cup_nU(J_n)$, but $S
\ssetminus U(J_n)$ is not empty for any fixed $n$, so there exists a
countable subset $S_0 \subset S$ not lying entirely in any of
$U(J_n)$. However, since $\kappa$ is uncountable and $J \subset
P_\kappa(S)$ is cofinal, we must have $S_0 \subset X$ for some $X
\in J$, and then $X \in J_n$ for some finite $n$.
\end{exa}

To construct $\kappa$-filtered partially ordered sets, it is
convenient to use the following modification of
Definition~\ref{filt.def}.

\begin{defn}\label{cons.filt.def}
A partially ordered set $J$ is {\em strongly $\kappa$-filtered} if
any map $S \to J$ admits a cone $S^> \to J$ that is a conservative
map. A conservative map $f:J \to J'$ between partially ordered sets is {\em
  weakly $\kappa$-filtered} if for any commutative diagram
\begin{equation}\label{adm.sq}
\begin{CD}
  S @>{e}>> S^>\\
  @V{a}VV @VV{g}V\\
  J @>{f}>> J'
\end{CD}
\end{equation}
with discrete $\kappa$-bounded $S$ and conservative $g$, there
exists a map $b:S^> \to J$ such that $b \circ e = a$ and $f \circ b
= g$.
\end{defn}

Note that with this definition, for any weakly $\kappa$-filtered map
$J \to J'$ with strongly $\kappa$-filtered target, $J$ is still
strongly $\kappa$-filtered. A $\kappa$-filtered set $J$ is strongly
$\kappa$-filtered iff it does not have maximal elements (since then
for any cone $S^> \to J$ with some vertex $j$, we have a non-trivial
map $f:j \to j'$, and $f_!S_>$ is conservative). On the other hand,
any surjective map $J \to \ppt$ is trivially weakly
$\kappa$-filtered (since there are no conservative maps $S^> \to
\ppt$ with non-empty $S$).

\begin{lemma}\label{sharp.le}
Assume given regular cardinals $\kappa \trl \mu$, and a
$\mu$-bounded set $S$. Then for any left-$\kappa$-bounded
$\mu$-bounded well-ordered set $Q_\kappa$, there exists a
$\mu$-bounded left $\kappa$-bounded partially ordered set $J$ and a
weakly $\kappa$-filtered conservative map $J \to Q_\kappa$ such that
the fiber $J_o$ over the smallest element $o \in Q_\kappa$ is
identified with the set $S$. Moreover, any functor $J_o \to I$ to a
$\kappa$-filtered category $I$ extends to a functor $J \to I$.
\end{lemma}

\proof{} Fix a set $S_\mu$ of cardinality $|S_\mu|=\mu$, and let $P$
be the set of pairs $\langle Q,J \rangle$ of a left-closed subset $Q
\subset Q_\kappa$, and a $\mu$-bounded subset $J \subset S_\mu
\times Q$ equipped with a partial order such that $J$ is
left-$\kappa$-bounded, and the projection $J \to Q$ is a
conservative order-preserving weakly $\kappa$-filtered map. In
particular, for any $\mu$-bounded subset $S \subset S_\mu$, $P$
contains the pair $\langle \{o\},S\rangle$, where $o \in Q$ is the
initial element, and $S$ carries the discrete order. Moreover, order
$P$ by setting $\langle Q',J' \rangle \leq \langle Q'',J''\rangle$
iff $Q' \subset Q''$, and $J' = (S \times Q') \cap J'' \subset S
\times Q''$, with the order induced from $J''$. Then since $\mu$ is
regular, $P$ with this order has upper bounds of all ascending
chains given by the union of the chain, so by the Zorn Lemma, it has
a maximal element $\langle J,Q \rangle \geq \langle \{o\},S
\rangle$. If $Q$ is not the whole $Q_\kappa$, then $Q = Q_\kappa/'q$
for some $q \in Q_\kappa$. Then consider $J$ with discrete order,
take the cofinal $\mu$-bounded subset $U \subset P_\kappa(J)$
provided by Definition~\ref{sharp.def}, forget the order on $U$,
choose an embedding $U \subset S_\mu$, and let $\lambda:U \to L(J)$
be the map sending some $\kappa$-bounded $X \subset J$ lying in $U
\subset P_\kappa(J)$ to the union $\lambda(X) = \cup_{x \in X}J/x
\subset J$. Then $Q'=Q_\kappa/q \cong Q^>$, and the union $J' = J
\copr U \subset (S \times Q) \copr S \cong S \times Q'$ with the
glued order $J' = \langle J,U,\lambda \rangle$ gives an element
$\langle Q',J' \rangle$ in $P$ strictly bigger than $\langle
J,Q\rangle$. This contradicts maximality.

If we are given a functor $J_o \to I$, then by the same induction,
it suffices to show that for any $q \in Q$, a functor $J/q
\ssetminus J_q \to I$ extends to $J/q$; this immediately follows
from the fact that $I$ is $\kappa$-filtered.
\endproof

\begin{prop}\label{sharp.prop}
Assume given regular cardinals $\kappa \trl \mu$. Then for any
$\mu$-bounded set $S$, there exists a $\mu$-bounded
left-$\kappa$-bounded $\kappa$-filtered partially ordered set $J_\mu$
equipped with a left-closed embedding $S \to J_\mu$ from $S$ equpped
with the discrete order. Moreover, any functor $S \to I$ to a
$\kappa$-filtered category $I$ extends to $J_\mu$.
\end{prop}

\proof{} Take the well-ordered left-$\kappa$-bounded set $Q_\kappa$
of Example~\ref{ord.exa}, construct $J$ and the conservative
weakly $\kappa$-filtered map $J \to Q_\kappa$ provided by
Lemma~\ref{sharp.le}, and note that since $Q_\kappa$ has no
maximal elements, it is strongly $\kappa$-filtered, so that
$J_\mu=J$ is $\kappa$-filtered.
\endproof

\begin{corr}\label{sharp.corr}
Assume given regular cardinals $\kappa \trl \mu$. Then any
left-$\kappa$-bounded $\mu$-bound\-ed $J \in \Pos_\mu$ admits a
left-closed embedding $J \subset J'$ into a left-$\kappa$-bounded
$\mu$-bounded $\kappa$-filtered $J' \in \Pos_\mu$. Moreover, any
functor $J \to I$ to a $\kappa$-filtered category $I$ extends to
$J'$.
\end{corr}

\proof{} Take a $\mu$-bounded cofinal subset $U \subset
P_\kappa(J)$, let $S$ be $U$ with discrete order, and take an
embedding $\eps:S \to J_\mu$ provided by
Proposition~\ref{sharp.prop}. Let $\lambda = \Lambda \circ \ups
\circ \eps^\dg:J_\mu \to L_\kappa(J)$, where $\ups:P_\kappa(S) \to
P_\kappa(J)$ sends a subset $S' \subset S \subset P_\kappa(J)$ to
the union of its elements, $\kappa$-bounded if so was $S'$, and
$\Lambda$ is the map \eqref{lp.eq} that sends $P_\kappa(J)$ into
$L_\kappa(J)$ since $J$ is left-$\kappa$-bounded. Then $J' = J
\copr_\lambda J_\mu$ does the job. Indeed, any $S' \subset J \subset
J'$ has an upper bound $j(S') \in S \subset J_\mu \subset J'$ since
$U \subset P_\kappa(J)$ is cofinal, and then for any $S' = S'_0
\copr S'_1 \in J' = J \copr_\lambda J_\mu$, take an upper bound for
$\{j(S_0')\} \cup S_1'$.
\endproof

\section{Ind-completions.}

Now assume given a category $I$, and again fix a regular cardinal
$\kappa$. Then by definition, the {\em $\Ind_\kappa$-completion}
$\Ind_\kappa(I)$ is the category of pairs $\langle J,i_\idot
\rangle$ of a $\kappa$-filtered small category $J$ and a functor
$i_\idot:J \to I$, with morphisms from $\langle J,i_\idot \rangle$
to $\langle J',i'_\idot \rangle$ given by
\begin{equation}\label{ind.hom}
\Hom(\langle J,i_\idot \rangle,\langle J',i'_\idot \rangle) =
\lim_{j \in J^o}\colim_{j' \in J'}\Hom(i_j,i'_{j'}).
\end{equation}
Objects of $\Ind_\kappa(I)$ are known as {\em $\Ind_\kappa$-objects}
in $I$, or simply $\Ind$-objects if $\kappa$ is the countable
cardinal. For any $I$, the category $\Ind_\kappa(I)$ is
$\kappa$-f-cocomplete. We have the fully faithful embedding $\eps:I
\to \Ind_\kappa(I)$, $i \mapsto \langle \ppt,i\rangle$, and it has a
left-adjoint $\eps_\dg:\Ind_\kappa(I) \to I$ if and only if $I$ is
$\kappa$-f-cocomplete, given by $\eps_\dg(\langle J,i_\idot \rangle)
= \colim_Ji_\idot$. In general, any functor $\gamma:I \to \C$ to a
$\kappa$-f-cocomplete $\C$ canonically extends to a functor $\gamma'
= \eps_\dg \circ \Ind_\kappa(\gamma):\Ind_\kappa(I) \to
\Ind_\kappa(\C) \to \C$ that preserves $\kappa$-filtered colimits,
and if $\gamma(i) \in \C$ is $\kappa$-compact for any $i \in I$, and
$\gamma$ is fully faithful, then \eqref{ind.hom} immediately shows
that $\gamma'$ is also fully faithful. In particular, if $I$ is
small, then $I^o\Sets$ is cocomplete, and the Yoneda embedding
\eqref{yo.eq} provides a full embedding
\begin{equation}\label{yo.filt.eq}
\Ind_\kappa(I) \to I^o\Sets.
\end{equation}
Here is what its essential image looks like.

\begin{lemma}\label{yo.ind.le}
For any essentially small category $I$, a functor $X:I^o \to \Sets$
lies in the essential image of the embedding \eqref{yo.filt.eq} iff
the category of elements $IX$ is $\kappa$-filtered.
\end{lemma}

\proof{} The ``if'' part immediately follows from
\eqref{yo.X.eq}. For the ``only if'' part, assume that $X =
\colim_J\Y(i_j)$ for some $\kappa$-filtered small $J$ and functor
$i_\idot:J \to I$. Then for any $j \in J$, we have a functor
$\eps_j:I\Y(i_j) = I / i_j \to IX$ sending the terminal object
$\id:i_j \to i_j$ to some $\langle i_j,x_j \rangle \in IX$, and for
any $E:IX \to \Sets$, we have
\begin{equation}\label{colim.IX}
\colim_{IX}E \cong \colim_J\colim_{I / i_j}\eps_j^*E \cong
\colim_JE(\langle i_j,x_j\rangle).
\end{equation}
Since $J$ is $\kappa$-filtered, the right-hand side of
\eqref{colim.IX} commutes with $\kappa$-bounded limits by
Proposition~\ref{filt.prop}, and by the same
Proposition~\ref{filt.prop}, this implies that $IX$ is
$\kappa$-filtered.
\endproof

\begin{corr}\label{small.le}
For any regular cardinals $\kappa' > \kappa$, a $\kappa$-bounded
Karoubi-closed small category $I$ is $\kappa'$-f-cocomplete.
\end{corr}

\proof{} It suffices to show that the Yoneda embedding $\Y:I \to
\Ind_{\kappa'}(I)$ is essentially surjective. Indeed, assume given
$X \in \Ind_{\kappa'}(I) \subset I^o\Sets$. Then for any $i \in I$,
$X(i)$ is $\kappa$-bounded by Corollary~\ref{sets.ka.corr}, so that
$IX$ is $\kappa$-bounded, and it is $\kappa$-filtered by
Lemma~\ref{yo.ind.le}. Therefore $\id:IX \to IX$ admits a
cone. Since $I$ is Karoubi-closed, $IX$ is also Karoubi-closed, and
then it has a terminal object by Lemma~\ref{P.le}, so that $X$ is in
the essential image of $\Y$ by Example~\ref{Y.exa}.
\endproof

\begin{remark}
Note that if a category $\C$ already has some $\kappa$-filtered
colimits, the embedding $\C \to \Ind_\kappa(\C)$ does not preserve
them --- effectively, we forget all the colimits we might have
already and formally add new ones. Thus unlike, say, Karoubi
closure, $\Ind_\kappa$ is not an idempotent operation, and it is
certainly not true that $\Ind_\kappa(\Ind_\kappa(\C)) \cong
\Ind_\kappa(\C)$.
\end{remark}

Any functor $\gamma:I' \to I$ between essentially small categories
induces a functor $\Ind_\kappa(\gamma):\Ind_\kappa(I') \to
\Ind_\kappa(I)$; in terms of the embedding \eqref{yo.filt.eq},
$\Ind_\kappa(\gamma)$ is given by the left Kan extension
$\gamma^o_!$ that obviously sends $\Ind_\kappa$-objects to
$\Ind_\kappa$-objects. The adjoint functor $\gamma^{o*}$ need not do
this. For a useful example of this phenomenon, assume given a
$\kappa$-bounded finite-dimensional partially ordered set $J$, and a
functor $\C:J \to \Cat$. Define $\Ind_\kappa(\C|J)$ as the full
subcategory $\Ind_\kappa(\C|J) \subset \Fun(\la{I\C^o},\Sets)$ whose
restriction to $\C(j)^o \subset \la{I\C^o}$ is in
$\Ind_\kappa(\C(j)) \subset \C(j)^o\Sets$ for any $j \in J$.

\begin{lemma}\label{lax.ind.le}
For any functor $\C:J \to \Cat$ from a $\kappa$-bounded
finite-dimen\-si\-onal partially ordered set $J$, we have
$\Ind_\kappa(\C|J) \cong \Ind_\kappa(\ra{\lim}_{J^o}\C)$.
\end{lemma}

\proof{} Consider the relative Yoneda embedding \eqref{yo.lax}, and
note that for any $s \in \ra{\lim}_{J^o}\C$ and $j \in J$, $\Y_J(s)$
restricts to $\Y(s(j))$ on $\C(j)^o \subset \la{I\C^o}$, so that
$\Y_J$ factors through $\Ind_\kappa(\C|J)$. Moreover, the twisted
arrow category $\Tw(J)$ is also a $\kappa$-bounded partially ordered
set, and then by \eqref{yo.lax.kan}, $\Y_J(s)$ is a $\kappa$-bounded
colimit of $\kappa$-compact objects, thus $\kappa$-compact by
Proposition~\ref{filt.prop}. Therefore $\Y_J$ extends to a fully
faithful functor $\Ind_\kappa(\ra{\lim}_{J^o}\C) \to
\Ind_\kappa(\C|J)$. By Lemma~\ref{yo.ind.le}, to see that it is
essentially surjective, we need to check that for any $X \in
\Ind_\kappa(\C|J) \subset \Fun(\la{I\C^o},\Sets)$, the
comma-category $\ra{\lim}_{J^o}\C \setminus_{\Y_J} X$ is
$\kappa$-filtered; this immediately follows from
Lemma~\ref{yo.ind.le} and Lemma~\ref{dim1.le}.
\endproof

\begin{remark}
Apart from the general works mentioned in the introduction, a useful
recent overview of $\Ind$-completions is \cite[Chapter 6]{KS}. In
particular, if $\kappa$ is the countable cardinal, and the functor
$\C:J \to \Cat$ is constant, Lemma~\ref{lax.ind.le} says that
$\Fun(J,-)$ commutes with $\Ind$; a slightly more general statement
can be found in \cite[Theorem 6.4.3]{KS}.
\end{remark}

\begin{corr}\label{lax.ind.corr}
For any small categories $\C_0$, $\C_1$, $\C$ equipped with functors
$\gamma_l:\C_l \to \C$, $l=0,1$, we have $\Ind_\kappa(\C_0
\vtimes_\C \C_1) \cong \Ind_\kappa(\C_0) \vtimes_{\Ind_\kappa(\C)}
\Ind_\kappa(\C_1)$.
\end{corr}

\proof{} Let the functor $\C_\idot:\V^o \to \Cat$ send $0$, $1$, $o$
to $\C_0$, $\C_1$, $\C$, with functors $\gamma_0$ and
$\gamma_1$. Then by Example~\ref{1.exa}, $\ra{\lim}_\V\C_\idot$ fits
into a cartesian square
\begin{equation}\label{lax.V.sq}
\begin{CD}
\ra{\lim}_\V\C_\idot @>>> \C_1 /_{\gamma_1} \C\\
@VVV @VV{\tau}V\\
\C_0 /_{\gamma_0} \C @>{\tau}>> \C,
\end{CD}
\end{equation}
where $\tau$ is as in \eqref{comma.facto}. We can rewrite
\eqref{lax.V.sq} as $\ra{\lim}_\V\C_\idot \cong \C_0 \vtimes_\C
\C'_1$, where $\C_1' = \C_1 /_{\gamma_1} \C$, and then $\eta$ of
\eqref{comma.facto} induces a fully faithful embedding
\begin{equation}\label{eta.v}
\id \vtimes_{\id} \eta:\C_0 \vtimes_\C \C_1 \to \C_0 \vtimes_\C
\C_1'.
\end{equation}
The same construction provides a fully faithful embedding
$$
\Ind_\kappa(\C_0) \vtimes_{\Ind_\kappa(\C)} \Ind_\kappa(\C_1) \to
\Ind_\kappa(\C_\idot|\V^o),
$$
and then the Yoneda embedding $\Y_{\V^o}$ of \eqref{yo.lax}
restricts to a fully faithful embedding $\C_0 \vtimes_C \C_1 \to
\Ind_\kappa(\C_0) \vtimes_{\Ind_\kappa(\C)} \Ind_\kappa(\C_1)$ whose
values are moreover $\kappa$-compact. Thus as in
Lemma~\ref{lax.ind.le}, it suffices to check that the induced fully
faithful functor $\Ind_\kappa(\C_0 \vtimes_\C \C_1) \to
\Ind_\kappa(\C_0) \vtimes_{\Ind_\kappa(\C)} \Ind_\kappa(\C_1)$ is
essentially surjective. As in Lemma~\ref{lax.ind.le}, this amounts
to showing that the corresponding comma-categories are
$\kappa$-filtered. Explicitly, for any $X = \langle X_0,X_1,\alpha
\rangle$ in $\Ind_\kappa(\C_0) \vtimes_{\Ind_\kappa(\C)}
\Ind_\kappa(\C_1) \subset \C_0^o\Sets \vtimes_{\C^o\Sets}
\C_1^o\Sets$, we have
$$
(\C_0 \vtimes_\C \C_1) /_\Y X \cong (\C_0X_0) \vtimes_{\C X} \C_1 X_1,
$$
where $X = \gamma^o_{1!}X_1:\C^o \to \Sets$, and then the functor
$\C_1 X_1 \to \C X$ is cofinal by Example~\ref{cof.exa}. We can now
replace $\C_0$, $\C_1$, $\C$ with $\C_0 X$, $\C_1 X_1$, $\C X$ to
simplify notation, and we just need to show that $\C_0 \vtimes_\C \C_1$
is $\kappa$-filtered as soon as $\gamma_1$ is cofinal. But then
$\eta:\C_1 \to \C_1'$ is cofinal by Corollary~\ref{comma.corr}, and
this immediately implies that the full embedding \eqref{eta.v} has
non-empty right comma-fibers. Since its target is $\kappa$-filtered
by Lemma~\ref{dim1.le}, we are done by Lemma~\ref{cof.le}.
\endproof

Let us now discuss what happens when one increases the
cardinality. For any $\kappa$-f-cocomplete category $\C$, denote by
$\Comp_\kappa(\C) \subset \C$ the full subcategory spanned by
$\kappa$-compact objects, and for any small category $\C$ and
regular cardinals $\kappa' \geq \kappa$, let
$\Ind_\kappa^{\kappa'}(\C) = \Comp_{\kappa'}(\Ind_\kappa(\C))$. Note
that for any regular cardinals $\mu \geq \kappa$, we have a fully
faithful embedding
\begin{equation}\label{mu.ka.eq}
\Ind_\mu(\Ind^\mu_\kappa(\C)) \to \Ind_\kappa(\C)
\end{equation}
induced by the full embedding $\Ind^\mu_\kappa(\C) \to
\Ind_\kappa(\C)$.

\begin{lemma}\label{ka.ka.le}
For any small category $\C$, $\Ind_\kappa^\kappa(\C) \cong P(\C)$ is
the Karoubi closure of the category $\C$.
\end{lemma}

\proof{} Since $\kappa$-f-cocomplete categories are Karoubi-closed,
and taking images of projectors commutes with colimits and colimits,
$P(\C) \subset \Ind_\kappa^\kappa(\C)$. Conversely, any
object $c \in \Ind_\kappa(\C)$ is by definition obtained as a
$\kappa$-filtered colimit $c = \colim_Ic_\idot$ for some $I$, and if
$c$ is $\kappa$-compact, $\id:c \to c$ must factor through $c_i$ for
some $i \in I$, so that $c$ is a retract of $c_i \in \C$.
\endproof

\begin{prop}\label{mu.ka.prop}
Assume given regular cardinals $\mu > \kappa$. Then \eqref{mu.ka.eq}
is an equivalence for any small category $\C$ if and only if $\kappa
\trl \mu$ in the sense of Definition~\ref{sharp.def}.
\end{prop}

\proof{} For the ``if'' part, by Lemma~\ref{yo.ind.le}, we need to
show that for any $X \in \Ind_\kappa(\C) \subset \C^o\Sets$, the
comma-category $\Ind_\kappa^\mu(\C) / X$ is $\mu$-filtered. We have
a full embedding $\Ind_\kappa^\mu(\C) \subset \Comp_\mu(\C^o\Sets)$,
and $\Comp_\mu(\C^o\Sets)/X$ is $\mu$-cocomplete by
Proposition~\ref{filt.prop}, thus $\mu$-filtered; therefore by
Lemma~\ref{cof.le}, it suffices to prove that for any $\mu$-compact
$Y \in \C^o\Sets$ equipped with a map $f:Y \to X$, the map factors
as $Y \to Y' \to X$ for some $Y' \in \Ind^\mu_\kappa(\C)$. Indeed,
by Corollary~\ref{filt.corr}, we may assume that $Y =
\gamma^o_!\ppt$ for some $\mu$-bounded left-finite partially ordered
set $J$ equipped with a functor $\gamma:J \to I$. The map $f$ then
defines a lifting $\gamma_f:J \to IX$ of the functor $\gamma$ to the
$\kappa$-filtered category $IX$, and by Corollary~\ref{sharp.corr},
there exists a $\mu$-bounded left-$\kappa$-bounded $\kappa$-filtered
partially ordered set $J'$ and a left-closed embedding $J \subset
J'$ such that $\gamma_f$ extends to a functor $\gamma_f':J' \to IX$,
or equivalently, to a functor $\gamma':J' \to I$ equipped with a map
$Y'={\gamma'}^o_!\ppt \to X$. This does the job.

For the ``only if'' part, note that for any set $S$, the partially
ordered set $P(S)$ of subsets $S' \subset S$ is cocomplete, with
colimits given by the union of subsets, and any $\kappa$-bounded $S'
\in P(S)$ is $\kappa$-compact. Now take a $\mu$-bounded set $S$,
consider the $\kappa$-filtered partially ordered set $J=P_\kappa(S)$
of Example~\ref{LP.exa}, and let $Q \subset P(J^>)$ be the set of
all $\kappa$-filtered subsets $X \subset J^>$, ordered by
inclusion. Note that for any $\kappa$-filtered $I$ and functor
$\gamma:I \to Q \subset P(J^>)$, any $\kappa$-bounded subset $Y \in
\colim_\gamma= \cup_{i \in I}\gamma(i)$, being $\kappa$-compact in
$P(J^>)$, must lie entirely within some $\gamma(i)$, so
$\colim_\gamma$ is $\kappa$-filtered. Therefore $Q$ is
$\kappa$-f-cocomplete.

Now, for any regular cardinal $\kappa'$, let $Q(\kappa') \subset Q$
be the subset of all $\kappa'$-bounded $\kappa$-filtered $X \subset
J^>$.  For any $X \subset J^>$, denote $X' = X \cup \{o\}$, and
observe that for any $X \in Q$ and cardinal $\kappa'$, we have $X'
\cong \colim_{P_{\kappa'}(X)}Y'$, where the limit is over all
$\kappa'$-bounded subsets $Y \subset X$. Therefore if $X$ is
$\kappa'$-compact, the inclusion $X \to X'$ must factor through some
$Y'$, so that $X \in Q(\kappa')$. Then in particular, $Q(\kappa)$
consists of $\kappa$-bounded $\kappa$-filtered subsets $X \subset
J^>$; these all must have a largest element, and conversely, any $X
\subset P_\kappa(J^>)$ that has a largest element is an element in
$Q(\kappa)$. For any $\kappa$-filtered $Y \subset J^>$, subsets $X
\in Q(\kappa)$ such that $X \subset Y$ are cofinal in $P_\kappa(Y)$,
so that $Y$ is a $\kappa$-filtered colimit of such subsets; we
conclude that $Q \cong \Ind_\kappa(Q(\kappa))$. But we then have
$\Ind^\mu_\kappa(Q(\kappa)) \subset Q(\mu)$, so by our assumption,
any $X \in Q$ is a $\mu$-filtered colimit of $\mu$-bounded
$\kappa$-filtered subsets $Y \subset J^>$.

In particular, $J \subset J^>$ is $\kappa$-filtered, thus lies in
$Q$. Therefore we must have $J = \colim_IY$ for some $\mu$-filtered
$I$ and functor $Y:I \to Q$, and since $o \in J^>$ is not in $J^>$,
all $Y(i) \in Q$ must also lie in $J \subset J^>$. Then for any $s
\in S$, $\{s\} \subset S$ is an element in $J$, so we must have
$\{s\} \in Y(i)$ for some $i \in I$, and then since $I$ is
$\mu$-filtered and $S$ is $\mu$-bounded, we must have some $Y=Y(i)
\in Q$ that contains all $\{s\}$. Since $Y$ is $\kappa$-filtered,
this means that it is cofinal in $J = P_\kappa(S)$.
\endproof

\begin{corr}\label{mu.ka.corr}
For any regular cardinals $\mu \geq \kappa$ and essentially small
category $\C$, $\Ind_\kappa^\mu(\C)$ is small.
\end{corr}

\proof{} If $\mu=\kappa$, the claim immediately follows from
Lemma~\ref{ka.ka.le} (the Karoubi closure of an essentially small
category is essentially small). If $\mu > \kappa$, then it suffices
to prove the claim for a cofinal class of cardinals $\mu$, and then
by Example~\ref{P.ka.exa}, we may assume that $\kappa \trl \mu$ in
the sense of Definition~\ref{sharp.def}. But then as we saw in the
proof of Proposition~\ref{mu.ka.prop}, any $\mu$-compact $X \in
\Ind_\kappa(\C)$ is a retract of an $\Ind_\kappa$-object $\langle
J,c_\idot \rangle \in \Ind_\kappa(\C)$ such that $J \in \Pos_\mu$ is
a $\mu$-bounded $\kappa$-filtered partially ordered set.
\endproof

\begin{remark}
Let us emphasize that even more so than in the rest of the paper,
the material in this Section and in Section~\ref{pos.sec} is due
entirely to \cite{AR}. This includes Definition~\ref{sharp.def} and
Proposition~\ref{mu.ka.prop}. Our proof of the ``only if'' part of
the latter is a simplication -- or some would say, bastardization --
of the proof of \cite[Theorem 2.11]{AR}. While the argument in
\cite{AR} is more general and more conceptual -- as a
counterexample, the authors use the whole category of
$\kappa$-filtered partially ordered sets rather than our $Q$ -- all
the essential points are borrowed from \cite{AR}, including the
wonderful idea of looking at $P(J^>)$ rather than the more obvious
$P(J)$.
\end{remark}

\section{Accessible categories.}

We are now ready to introduce our second main definition. Note that
for any regular cardinal $\kappa$ and $\kappa$-f-cocomplete category
$\C$, we have a fully faithful embedding
\begin{equation}\label{acc.eq}
\Ind_\kappa(\Comp_\kappa(\C)) \to \C
\end{equation}
induced by the fully faithful embedding $\Comp_\kappa(\C) \subset
\C$.

\begin{defn}\label{acc.def}
For any regular cardinal $\kappa$, a category $\C$ is {\em
  $\kappa$-accessible} if it is $\kappa$-f-cocomplete, the full
subcategory $\Comp_\kappa(\C) \subset \C$ of $\kappa$-compact
objects is essentially small, and the embedding \eqref{acc.eq} is an
equivalence. A functor $\gamma:\C \to \C'$ between
$\kappa$-accessible categories is {\em $\kappa$-accessible} if it
preserves $\kappa$-filtered colimits. A category $\C$ is {\em
  accessible} if it is $\kappa$-accessible for some regular cardinal
$\kappa$, and a functor $\gamma:\C \to \C'$ between accessible
categories is {\em accessible} if it is $\kappa$-accessible for some
$\kappa$ such that both $\C$ and $\C'$ are $\kappa$-accessible.
\end{defn}

Alternatively, a category $\C$ is $\kappa$-accessible iff $\C \cong
\Ind_\kappa(I)$ for some essentially small $I$ (and then by
Lemma~\ref{ka.ka.le}, $I$ is defined by $\C$ up to a Karoubi
closure). It even suffices to require that for some essentially
small $I \subset \Comp_\kappa(\C)$, the fully faithful functor
$\Ind_\kappa(I) \to \C$ is essentially surjective. By
Lemma~\ref{small.le}, any Karoubi-closed $\kappa$-bounded
essentially small category is $\kappa'$-accessible for any $\kappa'
> \kappa$. By Proposition~\ref{mu.ka.prop}, a $\kappa$-accesible
category is also $\mu$-accessible for any $\mu$ such that $\kappa
\trl \mu$, and then by Example~\ref{P.ka.exa}, the class of
cardinals $\kappa$ such that $\C$ is $\kappa$-accessible is cofinal;
moreover, since any $c \in \Ind_\kappa(I)$ is a small colimit, thus
$\mu$-compact for a large enough $\mu$, we actually have
\begin{equation}\label{comp.eq}
\C = \bigcup_\mu\Comp_\mu(\C),
\end{equation}
where the union is over all cardinals. All the categories
$\Comp_\mu(\C)$ in \eqref{comp.eq} are essentially small by
Corollary~\ref{mu.ka.corr}, and \eqref{comp.eq} implies that any
essentially small subcategory $\C' \subset \C$ lies inside
$\Comp_\mu(\C)$ for a large enough $\mu$. A $\kappa$-accessible
functor is also trivially $\mu$-accessible for any $\mu \geq
\kappa$, so that any collection of accessible categories and
functors is $\kappa$-accessible for a single cardinal $\kappa$ (and
then for a cofinal class of such cardinals).

\begin{lemma}
Assume given a pair of functors $\lambda:\C \to \C'$, $\rho:\C' \to
\C$ between accessible categories such that $\rho$ right-adjoint to
$\lambda$. Then both $\lambda$ and $\rho$ are accessible.
\end{lemma}

\proof{} The claim for $\lambda$ is trivial -- a functor that admits
a right-adjoint preserves all colimits, not only the filtered
ones. For $\rho$, choose $\kappa$ so that $\C$, $\C'$ and $\lambda$
are $\kappa$-accessible, and use \eqref{comp.eq} to find $\mu \geq
\kappa$ such that $\C$, $\C'$ and $\lambda$ are $\mu$-accessible,
and $\lambda$ sends $\Comp_\kappa(\C)$ into $\Comp_\mu(\C')$. Then
the full embedding $\C \cong \Ind_\kappa(\Comp_\kappa(\C)) \to
\Comp_\kappa(\C)^o\Sets$ is $\mu$-accessible, and its composition
with $\rho$ is $\mu$-accessible as well.
\endproof

By Corollary~\ref{lax.ind.corr}, for any regular cardinal $\kappa$,
$\Ind_\kappa$ commutes with lax fibered products, so that for any
$\kappa$-accessible categories $\C_0$, $\C_1$ and $\C$, and
$\kappa$-accessible functor $\gamma_l:\C_l \to \C$, $l=0,1$, the lax fiber
product $\C_0 \vtimes_\C \C_1$ is $\kappa$-accessible. The following
example shows that for the usual product $\C_0 \times_\C \C_1$, the
situation is more complicated: $\Ind_\kappa$ does not necessarily
preserve fibered products.

\begin{exa}\label{prod.exa}
Take $J = J_0 = J_1 = \N$, with the maps $q_l:\N \to \N$
sending $n$ to $2n+l$. Then we have the $\Ind$-object $\ppt \in
\Ind(\N)$, and both $\gamma_0$ and $\gamma_1$ are cofinal, so that
$\gamma^o_{l!}\ppt \cong \ppt$, $l=0,1$, and $\Ind(J_0)
\times_{\Ind(J)} \Ind(J_1)$ contains an object $\ppt \times
\ppt$. However, $J_0 \times_J J_1$ is empty.

More generally, for any regular $\kappa$, take a $\kappa$-filtered
left-$\kappa$-bounded well-ordered set $Q_\kappa$ of
Example~\ref{ord.exa}, and let $J_0=J_1=J=Q_\kappa \times \N$, with
lexicographical order. Then it is again a $\kappa$-filtered
left-$\kappa$-bounded well-ordered set, and we have cofinal
embeddings $(\id \times q_l):J_l \to J$, $l=0,1$, so
$\Ind_\kappa(J_0) \times_{\Ind_\kappa(J)} \Ind_\kappa(J_1)$ contains
an object $\ppt$, while $J_0 \times_J J_1$ is again empty.
\end{exa}

\begin{prop}\label{prod.prop}
For any regular cardinal $\kappa$ with successor cardinal
$\kappa^+$, $\kappa$-accessible categories $\C_0$, $\C_1$, $\C$, and
$\kappa$-accessible functors $\gamma_l:\C_l \to \C$, $l=0,1$, the
product $\C_0 \times_\C \C_1$ is $\kappa^+$-accessible.
\end{prop}

\proof{} By Example~\ref{suc.exa} and Proposition~\ref{mu.ka.prop},
the categories $\C_0$, $\C_1$, $\C$ are $\kappa^+$-accessible, and
then so are the functors $\gamma_0$, $\gamma_1$. Let
$I=\Comp_{\kappa^+}(\C)$, $I_l=\Comp_{\kappa^+}(\C_l)$, $l=0,1$,
with the induced functors $\gamma_l:I_l \to I$, $l=0,1$. These
categories are essentially small, and we might as well treat them as
small. Note that $\C$ is $\kappa$-f-cocomplete, and $I \subset \C$
is closed under $\kappa^+$-bounded $\kappa$-filtered colimits. Thus
for any $\kappa^+$-bounded $\kappa$-filtered $J \in \Pos_{\kappa^+}
\cap \Filt_\kappa$, any functor $i_\idot:J \to I$ has a colimit
$i=\colim_Ji_\idot \in I$.  Moreover, for any $X \in \C$,
$\Hom(-,X)$ preserves the corresponding limit $i=\lim_{J^o}i^o_\idot
\in I^o$, so that by \eqref{IX.eq}, any functor $J \to IX$ to the
category of elements $IX$ also has a colimit. The same holds for
$I_0$, $I_1$, the functors $\gamma_0$, $\gamma_1$ preserve
$\colim_J$, and for any $X_l \in \C_l \subset I_l^o\Sets$ with
$X=\gamma^o_{l!}X$. $l=0,1$, so do the cofinal functors
$\gamma_l:I_lX_l \to IX$, $l=0,1$ of Example~\ref{cof.exa}.

Now, $\Comp_{\kappa^+}(\C_0 \vtimes_\C \C_1) \cong I_0 \vtimes_I
I_1$ by Corollary~\ref{lax.ind.corr}, the full subcategory $\C_0
\times_\C \C_1 \subset \C_0 \vtimes_\C \C_1$ is closed under
$\kappa^+$-filtered colimits, thus $\kappa^+$-f-cocomplete and $I_0
\times_I I_1 = (\C_0 \times_\C \C_1) \cap I_0 \vtimes_I I_1 \subset
\Comp_{\kappa^+}(\C_0 \times_\C \C_1)$, so it suffices to prove that
for any $X = \langle X_0,X_1,\alpha \rangle$ in $\C_0 \times_\C \C_1
\subset I_0^o\Sets \times_{I^o\Sets} I_1^o\Sets$, the comma-category
$(I_0 \times_I I_1) / X$ is $\kappa^+$-filtered. Moreover, we have a
full embedding $(I_0 \times_I I_1) / X \to (I_0 \vtimes_I I_1) / X$
with $\kappa^+$-filtered target, so by Lemma~\ref{cof.le}, it
suffices to prove that this embedding has non-empty right
comma-fibers. As in the proof of Corollary~\ref{lax.ind.corr}, we
can now replace $I$, $I_0$, $I_1$ with $IX$, $I_0X_0$, $I_1X_1$, and
we are reduced to proving the following:
\begin{itemize}
\item for any $\kappa^+$-filtered categories $I$, $I_0$, $I_1$ that
  have colimits $\colim_J$ for any $J \in \Pos_{\kappa^+} \cap
  \Filt_\kappa$, and any cofinal functors $\gamma_l:I_l \to I$ that
  preserve these colimits, the full embedding $I_0 \times_I I_1
  \subset I_0 \vtimes_I I_1$ has non-empty right comma-fibers.
\end{itemize}
To prove this, assume given an object $\langle i'_0,i'_1,\alpha'
\rangle \in I_0 \times_I I_1$. Consider the partially ordered sets
$J \cong J_0 \cong J_1$ of Example~\ref{prod.exa}, with the cofinal
embeddings $q_l:J_l \to J$, $l=0,1$, and note that all three sets
are $\kappa$-filtered and $\kappa^+$-bounded. Denote by
$J_\idot:\V^o \to \Pos_{\kappa^+} \cap \Filt_\kappa$ the
corresponding functor. Take the initial segment $J' = \{o\} \times
[1] \subset J = Q_\kappa \times \N$ of the ordinal $J$, with $J'_l =
J' \cap J_l \cong \ppt$, $l=0,1$, and note that $\langle
i'_0,i'_1,\alpha' \rangle$ defines a functor $i'_\idot:J'_\idot \to
I_\idot$, where $I_\idot:\V^o \to \Filt_{\kappa^+}$ corresponds to
$I$ and $\gamma_l:I_l \to I$, $l=0,1$. Assume for a moment that we
can extend $i'_\idot$ to a functor $i_\idot:J_\idot \to I_\idot$,
and note that this finishes the proof: we can take $i^\dg_l =
\colim_{J_l}i_l$, $l=0,1$, and since the embeddings $q_0$, $q_1$ are
cofinal and $\gamma_0$, $\gamma_1$ preserve $\colim_J$, we have an
isomorphism $\alpha:\gamma_0(i^\dg_0) \cong
\colim_{J_0}\gamma_0(i_0) \cong \colim_Ji \cong
\colim_{J_1}\gamma_1(i_1) \cong \gamma_1(i^\dg_1)$ that defines
$\langle i^\dg_0,i^\dg_1,\alpha \rangle \in I_0 \times_I I_1 \subset
I_0 \vtimes_I I_1$ and a map $\langle i'_0,i'_1,\alpha \rangle \to
\langle i^\dg_0,i^\dg_1,\alpha\rangle$.

It remains to construct the extended functor $i_\idot$. To do it, by
the same induction as in Lemma~\ref{sharp.le}, it suffices to assume
that for some $j \in J$, with the corresponding functor $J_\idot /'
j:\V^o \to \Pos$ with values $J/'j$, $J_0/'j$, $J_1/'j$, we already
have a functor $i_\idot:J /'j \to I_\idot$, and we need to extend it
to $J_\idot/j$. Since set-theoretically, $J = J_0 \copr J_1$, we
have $j \in J_l$ for $l=0$ or $1$; the argument is the same in both
cases, so let us assume $l=0$. Then $J_1/'j = J_1/j$, so that only
extension happens at $J_0$ and $J$: we have $J_0/j \cong
(J_0/'j)^>$, $J/j \cong (J/'j)^>$, we have to find cones $i_{0>}$,
$i_>$ for $i_0:J_0/'j \to I_0$ and $i:J/'j \to I$, with some
vertices $i^\hash_0$, $i^\hash$, and an isomorphism
$\gamma_0(i^\hash_0) \cong i^\hash$. The procedure for doing this is
essentially the same as in Corollary~\ref{comma.corr}. First, find a
cone $i_{0>}$, with some vertex $i^\flat_0$. Next, consider the
corresponding functor $\gamma_0 \circ i_{0>} \copr i:J_0 / j
\copr_{J_0 /'j} J /' j \to I$, and find a cone for this functor,
with some vertex $i^\flat$. This gives a cone $i_>$ for $i$, with
some vertex $i^\flat$, and a map $\alpha:\gamma_0(i^\flat_0) \to
i^\flat$. Then take a map $\langle f_0,f\rangle:\langle
i_0^\flat,i^\flat,\alpha\rangle \to \langle
i^\hash_0,\gamma_0(i^\hash_0),\id \rangle$ provided by
Corollary~\ref{comma.corr}, with some $i^\hash_0 \in I_0$, and
replace $i_{0>}$ resp.\ $i_>$ with the cones $f_{0!}i_{0>}$
resp.\ $f_!i_>$.
\endproof

Finally, let us discuss what happens with the functor categories
(which, at the end of the day, was the whole point of the
exercise). For any $\kappa$-accessible categories $\C$ and $\E$, a
$\kappa$-accessible functor $E:\C \to \E$ is determined by its
restriction to the essentially small full subcategory
$\Comp_\kappa(\C) \subset \C$, uniquely up to a unique
isomorphism. Therefore for any two such functors $E$, $E'$, there is
at most a set of morphisms $E \to E'$, so that $\kappa$-accessible
functors $\C \to \E$ form a well-defined category
$\Fun_\kappa(\C,\E) \cong \Fun(\Comp_\kappa(\C),\E)$. Thus all
accessible functors $\C \to \E$ between accessible categories $\C$,
$\E$ form a well-defined category
\begin{equation}\label{fun.eq}
\Fun_\dg(\C,\E) = \bigcup_\kappa\Fun_\kappa(\C,\E),
\end{equation}
where the union is over all cardinals such that $\C$ and $\E$ are
$\kappa$-accessible. If $\E$ is accessible and $\C$ is essentially
small, then $\Fun(\C,\E) \cong \Fun(P(\C),\E)$, and the Karoubi
closure $P(\C)$ is accessible; in this case, the union in
\eqref{fun.eq} stabilizes at $\kappa^+$ for any $\kappa$ such that
$P(\C)$ is $\kappa$-bounded and $\E$ is $\kappa$-accessible. In
general, the whole category $\Fun_\dg(\C,\E)$ need not be
accessible; however, it is true for any individual term in the union
\eqref{fun.eq}.

\begin{lemma}\label{fun.le}
For any regular cardinal $\kappa$ and $\kappa$-accessible categories
$\C$, $\E$, the category $\Fun_\kappa(\C,\E)$ is accessible.
\end{lemma}

\proof{} Choose a small category $I$ equivalent to the essentially
small subcategory $\Comp_\kappa(\C) \subset \C$, so that
$\Fun_\kappa(\C,\E) \cong \Fun(I,\E)$. Then by \eqref{V.sq} and
Proposition~\ref{prod.prop}, it suffices to prove the claim when $I
= J$ is left-finite partially ordered set of finite dimension. By
Proposition~\ref{mu.ka.prop} and Example~\ref{P.ka.exa}, there
exists a regular cardinal $\mu$ such that $\E$ is $\mu$-accessible
and $J$ is $\mu$-bounded. We claim that $\Fun(J,\E)$ is
$\mu$-accessible. Indeed, by \eqref{V.eq}, any $E:J \to \E$ that
factors through $\Comp_\mu(\E) \subset \E$ is $\mu$-compact, so it
suffices to show that any $E:J \to \E$ is a $\mu$-filtered colimit
of functors $J \to \Comp_\mu(\E)$. If we now let $I =
\Comp_\mu(\E)$, then by Lemma~\ref{yo.ind.le}, $E$ defines a functor
$X:J \times I^o \to \Sets$ such that the category of elements
$IX(j)$ is $\mu$-filtered for any $j \in J$, and it suffices to
check that $\Fun(J,I) / E \cong \Sec(J^o,(J^o \times I)X)$ is
$\mu$-filtered. This is Lemma~\ref{dim1.le}.
\endproof

Summarizing, one can phrase the moral of the story as follows: a
category is ``tame'' if its Karoubi closure is accessible in the
sense of Definition~\ref{acc.def}, and all categories that appear
``in nature'' are tame. The necessity of taking Karoubi closures is
slightly unpleasant -- in particular, Karoubi closure does not
commute with standard operations such as fibered products or taking
functor categories, so one has to actually always remember it
explicitly; versions of Proposition~\ref{prod.prop} or
Lemma~\ref{fun.le} for tame categories are not true. But this is a
relatively small price to pay.

Let us finish the paper with a word of caution though: ``appears in
nature'' in the paragraph above should be understood in a rather
philosophical sense. Namely, of all the categorical operations
described in this Section, one is conspicuously missing --- that of
passing to the opposite category. There is a good reason for this:
the opposite $I^o$ of an accessible category $I$ is very rarely
accessible. The simplest example when it is not is $I = \Sets$. Thus
from the point of view presented here, simply inverting all arrows
in a category is not a ``natural'' operation, and the fact that it
gives something reasonable should be carefully justified in each
individual case.

\bigskip

{\noindent
Affiliations:
\begin{enumerate}
\renewcommand{\labelenumi}{\arabic{enumi}.}
\item Steklov Mathematics Institute (main affiliation).
\item National Research University Higher School of Economics.
\end{enumerate}}

{\noindent
{\em E-mail address\/}: {\tt kaledin@mi-ras.ru}
}

\end{document}